\documentclass{bull-l}

 \usepackage{graphicx, float,  hyperref, xcolor, tikz}

\newtheorem{theorem}{Theorem}[section]

\theoremstyle{definition}

\theoremstyle{remark}

\numberwithin{equation}{section}

\usetikzlibrary{shapes,arrows,shadows,calc}
 \usepackage{tkz-graph}
\tikzstyle{vertex}=[circle, draw, inner sep=0pt, minimum size=6pt]

\begin{document}

\title{Missing digits, and good approximations}


\author{Andrew Granville}
\address{D{\'e}partment  de Math{\'e}matiques et Statistique,   Universit{\'e} de Montr{\'e}al, CP 6128 succ Centre-Ville, Montr{\'e}al, QC  H3C 3J7, Canada.}
   \email{andrew.granville@umontreal.ca}  
\thanks{Thanks to Dimitris Koukoulopoulos, Sun-Kai Leung and Cihan Sabuncu for their comments on a draft of this article, and to James Maynard for sharing his graphics. The author is partially supported by NSERC of Canada, both by a Discovery Grant and by a CRC.}
 
\subjclass[2020]{Primary }

\begin{abstract}
James Maynard has taken the analytic number theory world by storm in the last decade, proving several important and surprising theorems, resolving questions that had seemed far out of reach.  He is perhaps best known for his work on small \cite{May1} and large \cite{May2} gaps between primes (which were discussed, hot off the press, in   \cite{Gra}).  In this article we will discuss two other Maynard breakthroughs:

 --- Mersenne numbers take the form $2^n-1$ and so appear as $111\dots 111$ in base 2, having no digit `$0$'. It is a famous conjecture that there are infinitely many such primes. More generally it was, until Maynard's work, an open question as to whether there are infinitely many primes that miss any given digit, in any given base. We will discuss Maynard's beautiful ideas that went into partly resolving this question \cite{May3}.

---  In 1926, Khinchin gave remarkable conditions for when real numbers can usually be ``well approximated'' by infinitely many rationals. However Khinchin's theorem regarded 1/2, 2/4, 3/6 as distinct rationals and so could not be easily modified to cope, say, with approximations by fractions with prime denominators.
In 1941 Duffin and Schaefer proposed an appropriate but significantly more general analogy involving approximation only by reduced fractions (which is much more useful).  We will discuss its recent resolution by Maynard together with Dimitris Koukoulopoulos \cite{KM}.
\end{abstract}

\maketitle

This year's Current Events Bulletin highlights the work of the 2022 Fields medallists. In James Maynard's case there are a surprising number of quite different breakthroughs that could be discussed.\footnote{In my forthcoming textbook about the distribution of primes, starting from the basics, about one-sixth of the book is dedicated to various  Maynard theorems. This, in one of the oldest and most venerable subjects of mathematics.} In my 2014 CEB lecture I described the work of Yitang Zhang \cite{Zh} on bounded gaps between primes and noted that a first-year postdoc, James Maynard, had taken a different, much simpler but related approach, to also get bounded gaps \cite{May1} (and a similar proof had been found, independently, by Terry Tao, and given on his blog). Versions of both Zhang's proof and the Maynard-Tao proof appear in my article  \cite{Gra}, where it is also announced that Maynard had within months made another spectacular breakthrough, this time on the largest known gaps between consecutive primes \cite{May2} (and a rather different proof \cite{FG1} had been found by Ford, Green, Konyagin and Tao, the two proofs combining to give an even better result \cite{FG2}).
It has been like this ever since with Maynard, many breakthrough results, some more suitable for a broad audience, some of primary importance for the technical improvements. Rather than attempt to summarize these all, I have selected two quite different topics, in  both of which Maynard proved spectacular breakthroughs on questions that had long been stuck.   

\part{Primes missing digits}
Most integers have many of each of the digits, 0 through 9, in their decimal expansion, so integers missing a given digit, or digits, are rare, making them hard to analyze.  For example, there are $3^k$ integers up to $10^k$ having only $7,8$ and $9$ in their decimal expansion as there are 3 possibilities for each of the $k$ digits in the expansion.\footnote{So there are about $x^\alpha$ integers up to $x$ having only $7,8$ and $9$ in their decimal expansion, where $\alpha:=\frac{\log 3}{\log 10}=0.4771\dots$} When we begin to explore we find the  primes 
\[
7, 79, 89, 97, 787, 797, 877, 887, 977, 997,\dots
\] having only the digits $7,8$ and $9$ in their decimal expansions. Are there infinitely many such primes? It seems likely given how many we have already found but this question, and questions like it, have long been wide open, researchers finding it difficult to find a method to plausibly attack such problems (as we will discuss below).  Indeed it was only  recently that researchers succeeded on the following related but seemingly less difficult problems:

--\ In 2010   Mauduit and Rivat \cite{MR} finally resolved Gelfond's problem
that the sum of the base-$q$ digits of prime numbers are equidistributed in arithmetic progressions, for all $q>2$. 

--\ In 2015  Bourgain \cite{Bo} showed that there are the expected number of primes with $k$ binary digits, for which $[ck]$ of those digits have preassigned values (and see Swaenepoel \cite{Sw} for base-$q$).

Maynard simplified and (in some aspects) sharpened the  tools used in these proofs but also added a perspective, and a technical confidence, that allowed him to surmount some of the established technical barriers.  Here we will sketch his proof giving an asymptotic for the number of primes up to large $x$, missing \emph{one} given digit in base $q$ (for $q$ sufficiently large), though his proof can be extended to counting 
the number of primes missing no more than  $\frac 15\, q^{2/5}$  base-$q$ digits  (again, for $q$ sufficiently large).  His proof works best if the allowed digits   lie in an interval, and in that case he was able to count the number of primes whose digits come from any sub-interval of $[0,q-1]$ of length $\gg q^{4/5}\log q$.

We begin by discussing where we should expect to find primes, and how many there are:

\section{Primes in arithmetic sequences}

We believe that an arithmetically natural set $\mathcal A$ of integers  contains infinitely many primes unless there is an obvious reason why not (like say, if $\mathcal A$ is the set of even integers, or the set of values of a reducible polynomial). Well known examples include, 

\begin{itemize}
\item $\mathcal A$ is the set of all integers;
\item $\mathcal A$ is the set of all integers in a given arithmetic progression (like $a \pmod q$ with $(a,q)=1$);
\item $\mathcal A=\{ p+2: p \text{ is prime}\}$, which is a way to ask for twin primes;
\item $\mathcal A=\{ n^2+1: n\in \mathbb Z\}$.
\end{itemize}
The first two questions are resolved and we even know  an asymptotic estimate for how many such primes there are up to a given $x$, while the second two questions are (wide) open.

\subsection{Guessing at the number of primes in $\mathcal A$}
The prime number theorem asserts that there are  $\sim \frac x{\log x}$ primes $\leq x$   (so roughly 1 in $\log x$ of the integers around $x$ are prime).\footnote{To prove such a result it helps to include a weight $\log p$ at each prime $p$ and prove instead that  $\sum_{ p \text{ prime}, p\leq x}   \log p\sim x$, since $x$ is a more natural function to work with than $\int_2^x \frac {dt}{\log t}$ (which is a more precise approximation than $\frac x{\log x}$). The prime number theorem can be deduced by the technique of ``partial summation'' which allows one to multiply or divide the summand by smooth weights.}
 As a first guess we might think that the primes are equidistributed amongst the arithmetic progressions mod $q$ and so the answer to the second question is $\sim \frac 1q \cdot \frac x{\log x}$; however $(a,q)$ divides any element of $a\pmod q$ and so if $(a,q)>1$ then this arithmetic progression contains at most one prime. Therefore we should restrict our attention to $a$ with $(a,q)=1$. There are $\phi(q)$ such progressions, and so we should adjust our guess so that if $(a,q)=1$ then there are $\sim \frac 1{\phi(q)} \cdot \frac x{\log x}$ primes $\leq x$ that are $\equiv a\pmod q$. This is the \emph{prime number theorem for arithmetic progressions}.\footnote{First claimed by de la Vall\'ee Poussin in 1899 based on ideas from his proof of the prime number theorem, and Dirichlet's proof of the infinitude of primes in arithmetic progressions. Thanks to Siegel and Walfisz this can be given, when $x$ is large enough compared to $q$, as follows: Fix reals $A,B>0$. If $q\leq (\log x)^A$ then the number of primes $\equiv a \pmod q$ up to $x$,
 \[
 \pi(x;q,a) = \frac{\pi(x)}{\phi(q)} \bigg( 1+ O\bigg( \frac 1{(\log x)^B}\bigg) \bigg)  \text{ whenever } (a,q)=1.
 \]}

Let $\mathcal A(x)$ be the set of integers in $\mathcal A$ up to $x$, and 
$\pi_{\mathcal A}(x)$ be the number of primes in $\mathcal A(x)$.
If the elements of $\mathcal A$ are as likely to be prime as random integers (roughly 1 in $\log x$  around $x$) then we'd guess that $\pi_{\mathcal A}(x)\approx  \frac{|\mathcal A(x)|}{\log x}$. This can be wrong since we have not accounted for any obvious biases in the set $\mathcal A$; for example, if $\mathcal A$ is the set of integers in an arithmetic progression mod $q$ then much depends on whether the progression is coprime with $q$. So we adjust our guess by a factor which is the probability that a random integer in $\mathcal A$ is coprime with $q$, divided by the probability, $\phi(q)/q$ that a random integer  is coprime with $q$. This then yields the guess
\[
\pi_{\mathcal A}(x) \approx \frac 1{\phi(q)/q}\cdot  \frac{|\mathcal A(x)|}{\log x} \sim \frac 1{\phi(q)} \cdot \frac x{\log x},
\]
when $(a,q)=1$ (and $0$ when $(a,q)>1$) so we recover the correct prediction.
This suggests a general strategy for guessing at $\pi_{\mathcal A}(x) $.

\subsection{Sparse sets of primes}
The first three questions above involve sets $\mathcal A$ that are quite dense amongst the integers. Our well-worn methods usually have limited traction with sets $\mathcal A$ that are \emph{sparse} such as
\begin{itemize}
\item $\mathcal A=\{ n\in (x,x+x^{.99}]\}$;
\item $\mathcal A=\{ n\equiv a \pmod q: n\leq x:=q^{100} \}$ for given  integer $q$ and $(a,q)=1$;
\item $\mathcal A=\{ n\leq x: \alpha n \pmod 1 \in [0,x^{-.01}] \}$ for a given real, irrational $\alpha$.
\end{itemize}
In each of these examples,  $|\mathcal A|\sim x^{.99}$, a rather sparse set. Each was shown to have more-or-less the expected number of primes  over 50 years ago (Theorems of Hoheisel, Linnik  and Vinogradov, respectively), though all known proofs are rather difficult. Moreover if we change ``$.99$'' to an exponent $<\frac 12$ then these questions are far beyond   our current state of knowledge.\footnote{The sparsest sets known in these questions to contain primes are $ (x,x+x^{.525}], x=q^5$ and $\alpha n \pmod 1 \in [0,x^{-\frac 13+\epsilon}]$ due to \cite{BHP, Xy, Mat} respectively.}

A family of sparse arithmetic sequences are given by the sets of values of polynomials (perhaps in several variables). Examples of sparse sets of values for which infinitely many primes have been found include

$\mathcal A= \{c^2+d^4: c,d\geq 1\}$ which has $|\mathcal A(x)|\asymp x^{3/4}$ (see \cite{FI}); and

$\mathcal A= \{a^3+2b^3: a,b\geq 1\}$ which has $|\mathcal A(x)|\asymp x^{2/3}$ (see \cite{HBL}).

\noindent This last set is an example of the set of values of a \emph{norm-form} as $a^3+2b^3$ is the norm an element, $a+2^{1/3}b$, of the ring of integers of $\mathbb Q(2^{1/3})$:

For a number field $K/\mathbb Q$, with ring of integers $\mathbb Z[\omega_1,\dots,\omega_d]$ the norm  
\[
\text{Norm}_{K/\mathbb Q}( x_1\omega_1+\ldots+x_d\omega_d) \in \mathbb Z[x_1,\dots,x_d]
\]
is a degree $d$ polynomial in the $d$ variables $x_1,\dots,x_d$. For example, $a+2^{1/3}b+2^{2/3}c$ has norm $a^3+2b^3+4c^3-6abc$.
The prime ideal theorem implies that norm-forms takes on infinitely many prime values with the $x_i$s all integers, provided that this polynomial  has no fixed prime factor.  These sequences are not so sparse (since they represent something like $x/(\log x)^C$ different integers up to $x$). 
However, in the last   example we know that there are roughly the expected number of prime values of  the norm-form for $a+2^{1/3}b+2^{2/3}c$  even when we fix $c=0$ (in which case one obtains a sparse set of integer values).

There are infinitely many prime values of the norm, $m^2+n^2$, of $m+in$ for   integers $m,n$, but if we fix $n=1$ we get the open question of primes of the form $m^2+1$. In 2002  Heath-Brown and Moroz \cite{HBM} proved that any cubic norm form with one of the variables equal to $0$ (as long as the new form is irreducible) takes roughly the expected number of prime values.  Moreover in 2018, Maynard  \cite{May4} proved such a result for norms of 
\[
\sum_{i=1}^{r} x_i\omega^i \in  \mathbb Z[\omega] \text{ where } [\mathbb Q(\omega):\mathbb Q]\leq \frac 43 r.
\]

Other than primes in short intervals, in short arithmetic progressions, and amongst polynomial values  perhaps the best known questions involving primes are those without some explicitly named digit or digits in their decimal or binary expansion:
\section{Primes with missing digits}
How many primes  only have the digits 1,\, 3 and 4 in their decimal expansions? When we start searching we find many:
\[
3, 11, 13, 31, 41, 43,   113, 131, 311, 313, 331,431,  433, 443, \dots
\]
and our  guess is that there are infinitely many such primes. To guess how many up to $x$, we can follow the above recipe: Here
 $|\mathcal A(10^k)|=3^k$,  and so
$|\mathcal A(x)|\asymp x^\alpha$ where $\alpha=\frac{\log 3}{\log 10}$.\footnote{Note that $|\mathcal A(3\cdot 10^k)|=|\mathcal A(2\cdot 10^k)|$ since no element of $\mathcal A$ begins with a `$2$', so  
$|\mathcal A(x)|/x^\alpha$ does not tend to a limit. The ratio ranges between  
$|\mathcal A(10^k)|/(10^k)^\alpha=1$ and $|\mathcal A(4\cdot 10^k)|/(4\cdot 10^k)^\alpha= 3 /4^\alpha=1.548\ldots $.}   We expect that the elements of $\mathcal A$ are independently equi-distributed modulo every prime $p$ except perhaps for those  dividing the base $10$:
Since the last digit of an element of $\mathcal A$ is $1,3$ or $4$, it is coprime with $10$ with probability $\frac 23$, whereas regular integers are coprime with $10$ with probability $\frac 12\cdot \frac 45=\frac 25$, and so
 we guess that   
 \[
 \pi_\mathcal A(x)\sim \frac {2/3}{2/5}   \cdot  \frac{|\mathcal A(x)|}{\log x} =\frac 53 \cdot  \frac{|\mathcal A(x)|}{\log x}.
 \]

\subsection{General prediction}
 If $\mathcal A$ is the set of integers $n$ which have only digits from $\mathcal D\subset \{ 0,1,\dots,q-1\}$ in their base $q$ expansion   let 
 $\mathcal D_q=\{ d\in \mathcal D: (d,q)=1\}$   and then we predict that 
  \[
 \pi_\mathcal A(x)\sim  \frac{ |\mathcal D_q| / |\mathcal D|}{\phi(q)/q} \cdot  \frac{|\mathcal A(x)|}{\log x},
 \]
 via the same reasoning.
Maynard proved this \cite{May3, May5} for certain general families of sparse sets $\mathcal A$.   His most spectacular result \cite{May3} yields (close to) the above with $q=10$ and   $|\mathcal D|=9$; that is, Maynard proved that there are roughly the expected number of primes that are missing \emph{one} given digit in decimal.\footnote{``Roughly'' meaning ``up to a multiplicative constant'' rather than an asymptotic.} His methods give a lot more (as we will describe). His methods can't handle sets as sparse as $\mathcal D=\{ 1,3,4\}$ with $q=10$; that is for another day.\footnote{Moreover there may be other, as yet undiscovered, reasons why there might not be any primes for a given $\mathcal D$. The one obstruction I know about is when every element of $\mathcal D$ is divisibly be some given prime $p$, which implies that all the elements of $\mathcal A$ are also divisible by $p$.}
  We will sketch slightly more than the easier argument from \cite{May5} which gives many results of this type though only for   bases that are significantly larger than $10$.

\subsection{Who cares?} 
Is this a silly question? It is certainly diverting to wonder whether there are infinitely primes with given missing digits, but how does that impact any other serious questions in mathematics?  This is a case of ``the proof of the pudding is in the eating'', that is, its real value   can be judged only from the beautiful mathematics that unfolds. The story is two-fold. The relevant Fourier coefficients have an extraordinary structure that allows Maynard to import ideas from Markov processes, and so prove such theorems in bases $>20000$. To get the base down to $10$, Maynard develops his ideas with a virtuosity in all sorts of deep techniques   that spin an extraordinary (though technical) tale.

\section{The circle method}

\subsection{Fourier analysis}
We use the identity
\[
\int_0^1 e(n\theta)d\theta = 1_{=0}(n) := \begin{cases}
1 & \text{ if } n=0;\\
0 & \text{ otherwise},
\end{cases}
\]
where $e(t):=e^{2 i \pi t}$ for any real $t$, and its discrete analog
\[
\frac 1N \sum_{j=0}^{N-1}  e\bigg( \frac {jn}N \bigg) =1_{=0}(n) \text{ whenever } |n|<N,
\]
obtained by summing the geometric series.

Let $\mathcal P$ denote the set of primes and $\mathcal A$  the set of integers missing some given digit or digits in  base-$q$.
To identify whether prime $p$ equals some $a\in \mathcal A$  we can take the above identities with $n=p-a$ and sum over all 
$a\in \mathcal A(N)$ and $p\in \mathcal P(N)$, to obtain, in the discrete case,
\begin{equation} \label{1}
\pi_{\mathcal A}(N) =  \sum_{p\leq N} \sum_{a\in {\mathcal A}(N)} \frac 1N \sum_{j=0}^{N-1}  e\bigg( \frac {j(p-a)}N \bigg) =
 \frac 1N \sum_{j=0}^{N-1} S_{\mathcal P}\bigg( \frac { j}N \bigg)S_{\mathcal A}\bigg( \frac {- j}N \bigg)
\end{equation}
where, for a given set of integers $T$, we define the \emph{exponential sum} (or the \emph{Fourier transform} of $T(N)$) by
\[
S_T(\theta) := \sum_{n\in T(N)}  e(n\theta) \text{ for any real } \theta. 
\]
Similarly, in the continuous case,  
\[
\pi_{\mathcal A}(N) =  \sum_{p\leq N} \sum_{a\in {\mathcal A}(N)} \int_0^1  e(  (p-a)\theta)d\theta =
 \int_0^1 S_{\mathcal P}(\theta)S_{\mathcal A}(-\theta)d\theta .
 \] 
One can work with either version, depending on whether discrete or continuous seems more convenient in the particular argument.\footnote{Some  of this discussion will make more sense to the novice if they think about the continuous version (though the discussion also applies to the discrete version).}
Writing 
\[
\pi_{\mathcal A}(N) =\sum_{n\leq N} 1_{\mathcal P}(n)1_{\mathcal A}(n)
\]
 one can also obtain \eqref{1}, and the continuous analog, from the Parseval-Plancherel identity in Fourier analysis.

\subsection{The circle method}
To establish a good estimate for $\pi_{\mathcal A}(N)$ using \eqref{1} one needs to identify those $j$ for which the summand on the right-hand side is large; for example, $S_T(0)=|T|$ and so the $j=0$ term in \eqref{1} yields 
\[ \frac 1N |\mathcal A(N)| \cdot \pi(N)\sim \frac{ |\mathcal A(N)| }{\log N }\]
 which is the expected order of magnitude of our main term (though it may be out by a multiplicative constant). Other terms where $\frac jN$ is small, or is close to a rational with small denominator often also contribute to the main term, whereas we hope that the combined contribution of all of the other terms is significantly smaller.  At first sight this seems unlikely since we only have the trivial bound $|S_T(\theta) |\leq |T|$ for the other terms, but the trick is to use the Cauchy-Schwarz inequality followed by Parseval's identity so that 
\[
\frac 1N \sum_{j=0}^{N-1}  | S_T( \tfrac { j}N)| \leq \bigg( \frac 1N \sum_{j=0}^{N-1}  |S_T( \tfrac { j}N)|^2 \bigg)^{1/2} =|T|^{1/2}.
\]
This implies for example that a typical term in the sum on the right-hand side of \eqref{1} has size $\sqrt{|{\mathcal A}(N)|}\cdot \sqrt{\pi(N)}$ which is  a little bigger than the main term but certainly not so egregiously as would happen if we used the trivial bound.

We have just described the basic thinking behind the \emph{circle method} used when one sums or integrates over the values of an exponential sum as the variable rotates around the unit circle (that is, $e ( \frac { j}N )$ for $0\leq j\leq N-1$, or $e(\theta)$ for $0\leq \theta<1$). When trying to estimate the sum on the right-hand side of \eqref{1}, we are most interested in those $\theta=\frac { j}N$ for which 
$S_{\mathcal P}(\theta)S_{\mathcal A}(-\theta)$ is ``large''.
Experience shows that with arithmetic problems, the exponential sums can typically only be large when $\theta$ is close to a rational with small denominators, and so we cut the circle up into these \emph{major arcs}, usually those $\theta$ near to a rational with small denominator, and 
 \emph{minor arcs}, the remaining $\theta$, bounding the contribution from the minor arcs, and being as precise as possible with the major arcs to obtain the main terms.
 
 Fourier analysis/the circle method is most successful when one has the product of at least three exponential sums to play with. For example the ternary Goldbach problem was more-or-less resolved by Vinogradov 85 years ago, whereas the binary Goldbach problem remains open.\footnote{It is known that \emph{almost all} integers $n$ can be written as the sum of two primes in the expected number of ways, since by counting over all integers $n$, one can estimate the variance via an integral involving three exponential sums.} 
 
 \subsection{The   ternary Goldbach problem} The number of representations of odd $N$ as the sum of three primes is given by 
 \[
 \int_0^1 e(-N\theta) S_{\mathcal P(N)}(\theta)^3 d\theta,
 \]
 and the arc of width $\asymp \frac 1N$ around $0$ yields a main term of size $\asymp \frac{N^2}{(\log N)^3 }$.
 We have the trivial bound $|S_{\mathcal P(N)}(\theta)|\leq \pi(N)$ and we will define here the minor arcs to be 
 \[
 \mathfrak m:= \{ \theta\in [0,1]: \ |S_{\mathcal P(N)}(\theta)|\leq \pi(N)/(\log N)^2\}.
 \]
 (Since the typical size of $|S_{\mathcal P(N)}(\theta)|$ is $\sqrt{\pi(N)}<N^{1/2}$ we expect that all but a tiny subset of the $\theta$ belong to these minor arcs.)
 Then 
 \begin{align*}
\bigg| \int_{\theta\in \mathfrak m} e(-N\theta) S_{\mathcal P(N)}(\theta)^3 d\theta \bigg| & 
\leq \int_{\theta\in \mathfrak m} |S_{\mathcal P(N)}(\theta)|^3 d\theta\\
&\leq  \frac{\pi(N)}{(\log N)^2} \cdot \int_{\theta\in [0,1)} |S_{\mathcal P(N)}(\theta)|^2 d\theta \\
&= \frac{\pi(N)^2}{(\log N)^2} \sim \frac{N^2}{(\log N)^4} 
 \end{align*}
 which is significantly smaller than the main term. Thus if we can identify which $\theta$ belong to $\mathfrak m$, then we can focus on evaluating $S_{\mathcal P(N)}(\theta)$ on the major arcs $\mathfrak M:=[0,1)\setminus \mathfrak m$. There are strong  bounds known for 
 $S_{\mathcal P(N)}(\theta)$, as we will see later, so these ambitions can all be achieved in practice.

\subsection{Major and minor arcs} The usual way to dissect the circle is to pick a parameter $1<M<N$ and recall that, by Dirichlet's Theorem (see the discussion in Part II), for every $\alpha\in [0,1]$ there exists a reduced fraction $r/s$ with $s\leq M$ for which 
\[
\bigg| \alpha-\frac rs\bigg| \leq \frac 1{sM}
\]
(and the right-hand side is $\leq 1/s^2$). Therefore we may cover $[0,1]$ (and so cover the circle, by mapping $t\to e(t)$) with the intervals (arcs),
\[
\bigcup_{s\leq M} \bigcup_{\substack{0\leq r\leq s\\ (r,s)=1}} \ \bigg[ \frac rs-\frac 1{sM},\ \frac rs+\frac 1{sM}\bigg].
\]
The arcs with $s$ small are usually the major arcs, those with  $s$ large are the minor arcs.

In our problem the partition of major and  minor arcs will be a bit more complicated.
The major arcs will be given by 
\[
\bigcup_{s\leq (\log N)^A} \bigcup_{\substack{0\leq r\leq s\\ (r,s)=1}} \     \bigg[ \frac rs-\frac {(\log N)^A}N ,\ \frac rs+\frac {(\log N)^A}N\bigg] ,
\]
and the main term will be obtained from those major arcs for which the prime factors of $s$ all divide $q$.
The minor arcs with be determined from the arcs above with $M=[\sqrt{N}]$, and then removing the major arcs.

Of course there is far more to say on the circle method than the brief discussion in this article. The reader should look into the two classic books on the subject \cite{Dav, Vau} for much more detail, and for applications to a wide variety of interesting questions.

  \section{The missing digit problem}
  
  Throughout let $\mathcal A$ be the set of integers whose digits come from the set $\mathcal D\subset \{ 0,1,\dots,q-1\}$. Our aim is to estimate $\pi_{\mathcal A}(N)$, and it will be convenient to let $N=q^k$ for some large even integer $k$.\footnote{For other large $N$ the key ideas are the same, but  dull technicalities arise.}
 
 The major arcs are typically given by the points $\theta\in [0,1)$ for which the integrand is large.\footnote{That is the goal, but one may have to include other points that one cannot easily exclude.}
If $S_{\mathcal P}(\theta)S_{\mathcal A}(-\theta)$ is large then $S_{\mathcal P}(\theta)$ and $S_{\mathcal A}(-\theta)$ must both individually be large. As we will see, Vinogradov proved that $S_{\mathcal P}(\theta)$ is only large when $\theta$ is near to a rational with small denominator.  $S_{\mathcal A}(\theta)$ behaves differently;   it is only large when there are many $0$'s and $q-1$'s in the decimal expansion of 
 $\theta$.  The simplest $\theta$ that satisfy both criteria take the form $\theta=\frac i{q^\ell}$ for some small $\ell$, perhaps with $\ell=1$ or, if $\ell>1$ then $\frac i{q^\ell}=\frac rs$, so that all the prime factors of $s$ must divide $q$.
   We therefore split the major arcs into three parts: Those   $\frac jN=\frac j{q^k}$ with
 \[
 \bigg|\frac jN-\frac rs\bigg|\leq \frac {(\log N)^A}N \text{ for some } 0\leq r\leq s \leq (\log N)^A \text{ with } (r,s)=1,
 \]
for some fixed $A>1$ where\smallskip
 
 ---    $s$ divides $q$, which contributes the main term;
 
 ---  $s$ only has prime factors which divide $q$ (excluding the  $\frac jN$ from the first case);
 
 --- $s$ is divisible by a prime not dividing $q$.
 \medskip
 
 \noindent We remark that $|\frac jN-\frac rs|\leq \frac {(\log N)^A}N$ if and only if  $|j- \frac rs\, N|\leq  (\log N)^A$.

 \subsection{The primary major arcs}

 Surprisingly the main term (in the discrete formulation) is obtained by simply taking those   $\theta=\frac j{q^k}$ for which $\theta= \frac \ell q$ for some integer $\ell$ (where $\ell$ and $q$ are not necessarily coprime). The contribution of such points to the above sum is
 \begin{align*}
q^{-k}  \sum_{\ell =0}^{q-1} S_{\mathcal P}\bigg( \frac { \ell }{q} \bigg)S_{\mathcal A}\bigg( \frac {-\ell }{q}\bigg) 
& =q^{-k}\sum_{a\in \mathcal A, a\leq q^k}     \sum_{ p \text{ prime}, \leq q^k}   \sum_{\ell =0}^{q-1} e\bigg(  \frac { \ell }{q} (p-a) \bigg)     \\
& =q^{1-k} \sum_{a\in \mathcal A, a\leq q^k}    \pi(q^k;q,a)  .
\end{align*}
Now if a prime $p$ does not divide $q$ and has last digit $d$ in base $q$ then $(d,q)=1$, and if $d\equiv p\equiv a \pmod q$ then $d\in \mathcal D$ so that 
$d\in \mathcal D_q$. There are $ |\mathcal D|^{k-1}$ integers $a\in \mathcal A, a\leq q^k$ with $a\equiv d \pmod q$, and so this sum becomes, using the prime number theorem for arithmetic progressions and as 
$|\mathcal A(q^k)|=|\mathcal D|^k$,
 \begin{align}
q^{1-k}  \sum_{d\in \mathcal D_q} |\mathcal D|^{k-1}   \pi(q^k;q,d) &\sim \frac{q^{1-k}  \cdot  |\mathcal A(q^k)| }{|\mathcal D|}
\sum_{d\in \mathcal D_q} \frac 1{\phi(q)} \frac{q^k}{\log q^k}  \notag\\ 
& =   \frac{ |\mathcal D_q| / |\mathcal D|}{\phi(q)/q} \cdot  \frac{|\mathcal A(N)| }{\log N},\label{eq: MainTerm}
\end{align}
which is precisely the prediction we had for $\pi_\mathcal A(N)$ above.

The asymptotic for $\pi_\mathcal A(N)$ now follows provided we can show that 
\begin{equation} \label{eq: Objective}
 \frac 1N \sum_{\substack{0\leq j\leq N-1 \\  \frac jN \ne \frac {r}q,\ 0\leq r\leq q-1}} \bigg| S_{\mathcal P}\bigg( \frac { j}N \bigg) \bigg|  \cdot  \bigg| S_{\mathcal A}\bigg( \frac {- j}N \bigg) \bigg| \ll
  \frac{|\mathcal A(N)| }{(\log N)^A}
\end{equation}
for some $A>1$. That is, we will be looking only at the absolute values of the exponential sums
$S_{\mathcal P}(\frac { j}N)$ and $S_{\mathcal A}( \frac {- j}N )$ and not trying to detect any surprising identities or cancelations based on angles.

 \subsection{Other major arcs, when all prime factors of $s$ divide $q$}  Throughout this subsection we assume that if prime $p$ divides $s$ then it divides $q$ so that $s$ divides $N=q^k$ for all sufficiently large $k$, and so $r/s$ may be written as $j/N$ for some integer $j$. We also assume that   $s\leq (\log N)^A$. 
 
 For these arcs we will find a strong upper bound on the values of $|S_{\mathcal P}( \frac jN)|$, and only bound $|S_{\mathcal A}( \frac jN)|\leq A(N)$, trivially:
 The prime number theorem for arithmetic progressions gives, if $(r,s)=1$,
 \begin{align}
 S_{\mathcal P}\bigg( \frac { r}{s} \bigg) &= \sum_{p\leq N} e\bigg( \frac { pr}{s} \bigg) = \sum_{a: (a,s)=1} e\bigg( \frac { ar}{s} \bigg)  \pi(N;s,a) +O(1)\notag  \\
 &  = \frac{\pi(N)}{\phi(s)} \sum_{a: (a,s)=1} e\bigg( \frac { ar}{s} \bigg)  + O\bigg( \frac {\pi(N)}{(\log N)^B}\bigg)
 \notag  \\ & =\pi(N) \bigg( \frac{\mu(s)}{\phi(s)}  + O\bigg( \frac {1}{(\log N)^B}\bigg) \bigg)  \label{eq: s-term}
\end{align}
as $\sum_{a: (a,s)=1} e(\frac { ar}{s})=\sum_{b: (b,s)=1} e( \frac {b}{s})=\mu(s)$ (an identity often credited to Ramanujan).
Therefore, by partial summation, if  $i$ is a non-zero integer with $  |i|\leq (\log N)^A$, or if $i=0$ and $\mu(s)=0$,
  \[
 S_{\mathcal P}\bigg( \frac { r}{s} +\frac { i}{q^k}\bigg)  = \pi(N)\frac{\mu(s)}{\phi(s)}  \int_0^N e \bigg( \frac{it}N\bigg) dt + O\bigg( \frac {i\pi(N)}{(\log N)^B}\bigg) \ll \frac {\pi(N)}{(\log N)^{B-A}}.
  \]
  We will write $\frac jN= \frac { r}{s} +\frac { i}{q^k}$ so that $ |i|\leq (\log N)^A$ if and only if $|j-\frac rs N|\leq (\log N)^A$.
 Therefore, since $|S_{\mathcal A}(\frac {- j}N)|\leq |\mathcal A(N)|$ trivially,  taking $B=4A-1$ with $A\geq 2$ we obtain
 \[
\frac 1N \sum_{\substack{s\leq (\log N)^A \\ p|s\implies p|q}} \sum_{\substack{0\leq r<s \\ (r,s)=1}} \sum_{\substack{ j: \\ \mu(s)^2\leq |j-\frac rs N|\leq (\log N)^A} }
\bigg|S_{\mathcal P}\bigg( \frac { j}N \bigg)S_{\mathcal A}\bigg( \frac {- j}N \bigg)\bigg| \ll \frac {|\mathcal A(N)|}{(\log N)^{A}}
 \]
 since there are $\ll (\log N)^A$ terms in the each of the sums. This upper bound
   is much smaller than the main term in \eqref{eq: MainTerm}.
 
 The only $r/s$ which are not accounted for here are those   where $s$ is squarefree and all  of its prime factors divide $q$. But this implies   that  $s$ divides $q$, and these terms were already  included in the sum in the previous subsection that led to \eqref{eq: MainTerm}.
 Therefore the calculations in this and the previous subsection account for the contributions to the sum in \eqref{1} of the ``$q$-smooth'' major arcs
 \[
\bigcup_{\substack{s\leq (\log N)^A \\  p|s \implies p|q}} \bigcup_{\substack{0\leq r\leq s\\ (r,s)=1}} \     \bigg[ \frac rs-\frac {(\log N)^A}N ,\ \frac rs+\frac {(\log N)^A}N\bigg].
\]
Before finishing with the major arcs we will need to introduce a key perspective for working with the exponential sums 
$|S_{\mathcal A}(\alpha)|$.

 \section{What makes restricted digit problems tractable?}
 From Parseval we know that for a given set $T$, we typically have $|S_T(\alpha)|\ll T(N)^{1/2}$
 (and for most $T$, we expect that $|S_T(\alpha)|\asymp T(N)^{1/2}$ for almost all $\alpha$).
 Therefore using Parseval  we have, for most $\alpha$, 
 \[
 |S_{\mathcal A}(\alpha)|\cdot |S_{\mathcal P}(\alpha)| \ll (A(N)\cdot \pi(N))^{1/2} \asymp N^{1-\delta+o(1)}
 \]
 where we define $\delta>0$ by $|\mathcal D|=q^{1-2\delta}$. However this is \emph{much bigger} than the main term $ \frac{|\mathcal A(N)| }{\log N}=N^{1-2\delta+o(1)}$, and so the circle method approach to digit sum problems has long seemed hopeless, since the sum of the absolute values of the contributions from the minor arcs seems likely to be so much larger than the main terms.

However, Maynard observed that the values of $|S_{\mathcal A}(\alpha)|$ are quite unusual in that they are not typically of size $A(N)^{1/2}$ but rather they are usually much  smaller, as we shall see.
 Therefore  restricted digit problems in base $q$ are more tractable because the structure of $\mathcal A$ leads to an unusual distribution of the sizes of its corresponding exponential sums, and so the contributions from the minor arcs are typically surprisingly small.

\subsection{The extraordinary structure of these exponential sums}

If $\mathcal A$ is the set of integers, whose base-$q$ digits come only from the set $\mathcal D\subset\{ 0,1,\dots, q-1\}$, and $N=q^k$, then we can write
 \[
 \mathcal A(N)=\bigg\{ n=\sum_{i=0}^{k-1} a_iq^i: \text{ Each } a_i\in \mathcal D\bigg\}.
 \]
 Since $e(n\theta) = \prod_{i=0}^{k-1} e(a_iq^i\theta)$, therefore
 \begin{align}
 S_{\mathcal A}(\theta) &= \sum_{\text{ Each } a_i\in \mathcal D }  \prod_{i=0}^{k-1} e(a_iq^i\theta)
 = \prod_{i=0}^{k-1}\bigg( \sum_{ a_i\in \mathcal D }   e(a_iq^i\theta)  \bigg) \notag\\
 & =  \prod_{i=0}^{k-1}\bigg( \frac{e(q^{i+1}\theta) -1} {e(q^i\theta) -1}  - e(b q^i\theta) \bigg)  \label{2}
 \end{align}
 where we have assumed that $\mathcal D=\{ 0,1,\dots, q-1\}\setminus \{b\}$ only in the last displayed line.
 It is very unusual for an exponential sum of interest to be a product of much simpler exponential sums like this. If the exponential sums in the product were independent of each other then we could focus on each $i$ separately and get best possible results; however the value of $q^i\theta \mod 1$ can be used to determine $q^{i+1}\theta \mod 1$ and so these are not independent. However, in practice, especially if $q$ is large, they will be independent enough to get some surprisingly strong upper bounds on $|S_{\mathcal A}(\theta)|$ for most $\theta$.

We define
\[
F_{\mathcal D}(\phi):=  \bigg| \sum_{n\in \mathcal D} e(n\phi)\bigg|\leq | \mathcal D|,
\]
so that
\[
|S_{\mathcal A}(\theta)|= | \mathcal A(N)| \cdot  \prod_{i=0}^{k-1} \frac 1{ | \mathcal D|} \, F_{\mathcal D}(q^i\theta)
\]
since $| \mathcal A(N)|=| \mathcal D|^k$.

\subsection{First upper bounds when $\mathcal D=\{ 0,1,\dots, q-1\}\setminus \{b\}$}
  Taking absolute values and using the triangle inequality we have 
  \begin{align} 
  F_{\mathcal D}(\phi)&=\bigg| \frac{e(q\phi) -1} {e(\phi) -1}  - e(b \phi) \bigg| \leq    1+ \frac{|e(q\phi) -1|} {|e(\phi) -1|}\notag \\
  &   \leq    1+ \frac{2} {|e(\phi) -1|}  = 1+ \frac 1{ \sin(\pi \| \phi \|)} \label{eq: sinbound}
   \end{align}
where $\| t\| =\min_{n\in \mathbb Z} |t-n|$, and therefore
  \begin{equation} \label{eq: absbound}
  F_{\mathcal D}(\phi) \leq 1+ \frac 1{ 2\| \phi \|}
  \end{equation}
 since   $  \sin (\pi \| t\|)\geq 2\| t\|$.
     
 Now if  $\theta=\sum_{j\geq 1}\frac{ t_{j-1}}{q^j}$  in base $q$ (with the $t_i\in \{ 0,1,\dots,q-1\}$) then
   \[
 q^i\theta \mod 1 = \frac{t_i}q+\frac{t_{i+1}}{q^2}   +\dots =\frac{ t_i+   (q^{i+1}\theta \mod 1)}q,
 \]
 and so $q^i\theta \mod 1\in [ \frac{t_i}q,\frac{t_i+1}q)$. This implies that 
    $\| q^i\theta \|\geq \min\{  \frac{t_i}q,1-\frac{t_i+1}q\}$ and so, by \eqref{eq: sinbound},
    \[
 F_{\mathcal D}(q^i\theta) \leq
   \min\bigg\{ q-1, 1+ \frac 1{  \min\{ \sin(\pi \frac {t_i}q), \sin(\pi \frac {q-1-t_i}q\} } \bigg\} ,
  \]
 and we obtain, in \eqref{2},
 \begin{equation} \label{3}
 | S_{\mathcal A}(\theta)| \leq   \prod_{i=0}^{k-1}  \min\bigg\{ q-1, 1+ \frac 1{  \min\{ \sin(\pi \frac {t_i}q), \sin(\pi \frac {q-1-t_i}q\} } \bigg\}.
 \end{equation}
In particular if, as is typical,  $q^{2/3}<t_i < q-q^{2/3}$ then the $i$th term in \eqref{2} is $\ll q^{1/3}$, a big improvement over the Parseval bound $\sqrt{q-1}$.

In fact for almost all $\theta$ the $t_i$ are uniformly distributed in $[0,q-1]$, that is 
$\# \{ i\in [1,k]: t_i=r\} \sim k/q$ for all $r\in [0,q-1]$, and so 
\[
|S_{\mathcal A}(\theta)|\ll \bigg( q \prod_{1\leq r\leq q/2} \bigg( 1+ \frac 1{ \sin(\pi \frac rq)} \bigg)\bigg)^{\{ 2+o(1)\} k/q}=(C+o(1))^k
\]
where $C:=\exp(\frac 4\pi L(2,(\frac{-4}\cdot))\approx 3.209912300$.
This is much smaller than $q^{k/2}$ for large $k$. As promised we have shown that the $|S_{\mathcal A}(\theta)|$, where $\mathcal A$ is the set of integers missing one particular digit in base $q$,  have a very different distribution from the $|S_T(\theta)|$ for a typical set of integers $T$.
This distribution indeed implies that the set of $\theta$ for which 
$|S_{\mathcal A}(\theta)|$ is not very small, has tiny measure.  We follow Maynard's argument to exploit this.
 
\subsection{Major arcs, where $s$ has a prime factor that does not divide $q$.}   \label{sec 5.3}

 A  weaker bound on the $i$th term, but which is easier to work with, comes from noting that 
 \[
 |e(a\phi) + e((a+1)\phi)|^2 = 2 + 2 \cos(2\pi \phi) < 4 \exp( - 2\| \phi\|^2),
 \]
 so that $ |e(a\phi) + e((a+1)\phi)|\leq 2 \exp( - \| \phi\|^2)$.
 If $q>3$ then there are two consecutive integers in $\mathcal D$ and so
 \[
 \sum_{ a\in \mathcal D }   e(a\phi) \leq q-3 + 2 \exp( - \| \phi\|^2) \leq (q-1) \exp\bigg( -\frac{\| \phi\|^2}q\bigg),
 \] 
 and therefore, by \eqref{2},
 \begin{equation} \label{eq: thetasum}
  | S_{\mathcal A}(\theta)| \leq |\mathcal A(N)| \exp\bigg( -\frac1q \sum_{i=0}^{k-1} \| q^i \theta\|^2\bigg)
 \end{equation}
 We use this not very good upper bound to deal with the (few) remaining possible major arcs, though these arguments, and so  \eqref{eq: thetasum}, can easily be sharpened.
 
 
 Suppose that prime $p|s$ but  $p\not |q$.  Then $p$ divides the denominator of  the reduced fraction for $q^i\cdot\frac rs$  so that $\| q^i\cdot \frac rs\| \geq \frac 1p$. Moreover if  $|\theta-\frac rs|\leq \frac 1{2pN^{1/2}}$ and $i\leq \frac k2$  then 
\[
\| q^i \theta\| \geq \| q^i\cdot \tfrac rs\|  - q^i |\theta-\tfrac rs| \geq \tfrac 1p - \tfrac {q^{k/2}}{2pN^{1/2}} = \tfrac 1{2p}.
 \]
 Now if $\| q^i \theta\| <\frac 1{2q}$ then $\| q^{i+1} \theta\| = q \| q^i \theta\|$. Therefore, for every integer $i$ there exists an integer $j, i\leq j\leq i+\lfloor \frac{\log p}{\log q}\rfloor $ for which  $\| q^j \theta\| \geq \frac 1{2q}$, which implies that 
 \[
 \sum_{i=0}^{k-1} \| q^i \theta\|^2 \geq   \sum_{i=0}^{k/2} \| q^i \theta\|^2 \geq \frac 1{4q^2} \#\{ j\in [0,\tfrac k2): \| q^j \theta\| \geq \tfrac 1{2q}\} \geq  \frac 1{4q^2}  \frac{\log q^{k/2}}{ \log pq } \geq  \frac k{8mq^2}  
 \]
 for $s\leq q^{m}$ and $m\in \mathbb Z$, since then $\lfloor \frac{\log p}{\log q}\rfloor \leq m-1$. 
 Here we let $m=\lfloor \sqrt{k}/9q^3 \rfloor$ and assume that $k\geq 100q^6$.
 
 Thus 
 $ | S_{\mathcal A}(\theta)| \leq |\mathcal A(N)| \exp ( - \frac k{8mq^3}) $ by    \eqref{eq: thetasum},
and $|S_{\mathcal P}(\theta)|\leq \pi(N)$ trivially, so that as $2q^{2m}\leq N^{1/2}$ then
\begin{align*}
\frac 1N \sum_{\substack{s\leq q^{m} \\ \exists p|s,  p\not |q}} \sum_{\substack{0\leq r<s \\ (r,s)=1}} \sum_{j:  |j-\frac rs N|\leq   q^{m}} 
\bigg|S_{\mathcal P}\bigg( \frac { j}N \bigg)S_{\mathcal A}\bigg( \frac {- j}N \bigg)\bigg| 
&\ll \frac{  |\mathcal A(N)| }{\log N} q^{3m}  \exp ( - \tfrac k {8mq^3}) \\ &\ll \frac{  |\mathcal A(N)| }{\log N} e^{-\sqrt{k}},
 \end{align*}
 which is much smaller than the main term in \eqref{eq: MainTerm}.

 This subsection accounts for the major arcs,
 \[
\bigcup_{\substack{s\leq q^m \\  \exists p|s \text{ such that } p\not| q}} \bigcup_{\substack{0\leq r\leq s\\ (r,s)=1}} \     \bigg[ \frac rs-\frac {q^m}N ,\ \frac rs+\frac {q^m}N\bigg].
\]
where $q^m=c_q^{\sqrt{k}}$ for some $c_q>1$, which is much  larger than $(\log N)^A$ for $k$  sufficiently large.

\section{The remaining challenge; the minor arcs}  \label{sec: Minor}

When dealing with each of the second two types of major arcs we bounded one of the exponential sums trivially; we will have no such luxury when bounding the contribution of the minor arcs. We obtain the minor arcs $\mathfrak m$, for $M=\lfloor \sqrt{N} \rfloor=q^{k/2}$, from subtracting the major arcs from a partition of the unit circle:
\[
  \bigcup_{\substack{0\leq r\leq s\leq M\\ (r,s)=1}} \ \bigg[ \frac rs-\frac 1{sM},\ \frac rs+\frac 1{sM}\bigg] \ \setminus \ 
  \bigcup_{\substack{0\leq r\leq s\leq (\log N)^A \\ (r,s)=1}} \     \bigg[ \frac rs-\frac {(\log N)^A}N ,\ \frac rs+\frac {(\log N)^A}N\bigg].
\]
We can further partition these   arcs according to the sizes of $s$:\footnote{$x\asymp X$ means $x$ runs through 
the integers or reals (as  appropriate) in the interval $(X,qX]$.}
\[
s\asymp S   \text{ with } 1\leq S=q^i \leq M/q
\]
where $i\geq 0$ is an integer, with $i\leq k/2$ (where $k$ is even); and   the size of  $\| s\theta\|$ for $\theta= \frac jN$:
\[
 \bigg| \frac jN -\frac rs\bigg|\leq \frac 1N,\text{ or }  \bigg| \frac jN -\frac rs\bigg|\asymp \frac BN  \text{ with } 1\leq B=q^\ell
\]
 where $\ell\geq 0$ is an integer, and so that 
\[
B=q^\ell\leq \frac N{q^2SM}  
\]
since $|\frac jN -\frac rs|< \frac 1{sM}$; that is, $i+\ell\leq \tfrac k2-2$.
 This also implies that
$\| s \frac jN \| = s \|  \frac jN  \| \leq \frac 1M$.

The major arcs took account of the cases in which both $B, S \ll (\log N)^A$, and so for the minor arcs we have
$B S \gg (\log N)^A$, so that
\[
(\log N)^A \ll BS \leq \frac N{q^2M}  
\]
(that is, $\log k\ll_q i+\ell \leq \frac k2-2$).

\subsection{Well-known bounds on $S_{\mathcal P}(\theta)$}
 Vinogradov's estimate for exponential sums (\cite{Dav1}, pg 142) 
 gives that if $\alpha=\frac jN=\frac rs + \beta$ with $(r,s)=1$ and $|\beta|<\frac 1{s^2}$ then 
 \[
 |S_{\mathcal P}(\alpha)| \ll \bigg(N^{4/5} + (s N)^{1/2}+ \frac{N}{ s^{1/2}}  \bigg) (\log N)^4\ll \bigg(N^{4/5} +   \frac{N}{ S^{1/2}}  \bigg) (\log N)^4
 \]
 since $(s N)^{1/2}\leq (M N)^{1/2}\leq N^{4/5}$ as $M\leq N^{3/5}$ and as $s\asymp S$. We use this in the first range above.
 
 In the second range above we have $\| s \frac jN\| \asymp  \frac{BS}N$.
 By a slight modification of Vinogradov's proof, we also have the bound
 \begin{align}
 |S_{\mathcal P}(\alpha)|& \ll \bigg(N^{4/5} + \frac{N^{1/2}}{ \|s \alpha\|^{1/2}} + \|s \alpha\|^{1/2} N\bigg) (\log N)^4 \notag\\
 & \ll \bigg( N^{4/5} +\frac N {(BS)^{1/2}}  \bigg) (\log N)^{4} \label{eq: ExpSumPrime}
\end{align}
since $ \|s \frac jN \|^{1/2} N\asymp (BSN)^{1/2}\ll \frac N{M^{1/2}} \leq N^{4/5}$ as $M\geq N^{2/5}$. 
 
\subsection{The end-game}
Our main goal in this section is to show that if $q\geq 133359$ and
$\mathcal D=\{ 0,1,\dots, q-1\}\setminus \{ b\}$
 then
\begin{equation} \label{KeyMinors}
 \sum_{\substack{0\leq r<s\leq S \\ (r,s)=1}}  \sum_{j: |j-q^k \cdot \frac rs|\leq B} \bigg| S_{\mathcal A}  \bigg( \frac j {q^k} \bigg)\bigg|      
\ll_q |\mathcal A(N)|      (BS^2)^{\frac 15-\eta},
 \end{equation}
for some   $\eta>0$, where the ``$\ll$'' depends only on $q$.\footnote{In this case, $A(N)=N^{1-\delta_q}$ for $N=q^k$ where $\delta_q=\frac{\log (1+\frac 1{q-1})}{\log q}$, so that the bigger that $q$ gets, the more (Hausdorff)-dense $\mathcal A$ is. This is why these arguments work better as $q$ gets larger.}
 Using the bound in \eqref{eq: ExpSumPrime} we then deduce that 
\[
 \sum_{\substack{0\leq r<s\asymp S \\ (r,s)=1}}  \sum_{\substack{|\frac jN- \frac rs|\leq \frac 1N \text{ or}\\ |\frac jN- \frac rs|\asymp \frac BN}} 
\bigg| S_{\mathcal P}  \bigg( \frac j N \bigg)  \cdot   S_{\mathcal A}  \bigg( \frac {-j} N \bigg)\bigg| 
\ll_q |\mathcal A(N)|   \bigg( N^{1-\eta} +\frac N {(BS)^{\frac 1{10}}}  \bigg) (\log N)^{4}   
\]
since $BS^2\leq BSM \ll N$ and $BS^2\leq(BS)^2$. Now we sum this bound over all $B=q^\ell, S=q^i$ where $i$ and $\ell$ are integers $\geq 0$, with $(\log N)^A \ll BS=q^{i+\ell} \leq N$ (so that there are $\ll (\log N)^2$ such pairs $i,\ell$).
Therefore we obtain
\[
\frac 1N  \sum_{j: \frac jN\in \mathfrak m} 
\bigg| S_{\mathcal P}  \bigg( \frac j N \bigg)  \cdot   S_{\mathcal A}  \bigg( \frac {-j} N \bigg)\bigg| 
\ll_q \frac{|\mathcal A(N)|}{(\log N)^C}
\]
provided $A\geq 10(C+4)$. (We can therefore define our arcs using any fixed $A>50$, and then select $C$ with 
$A= 10(C+4)$, ensuring that $C>1$.)
Therefore
Maynard's result, that we have asymptotically the predicted number of primes missing some given digit in base $q$,  follows for all bases $q\geq 133359$.

\subsection{The mean value of $|S_{\mathcal A}  (\alpha)|$}

For any real $\theta$ the set of values of the first $k$ base-$q$ digits of 
\[
\bigg\{ \theta +\frac j{q^k} \mod 1: 0\leq j\leq q^k-1\bigg\}
\]
 run once through each    $(t_0,\dots,t_{k-1})\in \{ 0,1,\dots, q-1\}^k$.  Therefore, by \eqref{3},
\begin{equation} \label{eq: s-sum}
\sum_{j=0}^{q^k-1}  \bigg| S_{\mathcal A} \bigg(\theta +\frac j{q^k} \bigg)  \bigg| \leq
\prod_{i=0}^{k-1} \sum_{t_i=0}^{q-1} \min\bigg\{ q-1, 1+ \frac 1{  \min\{ \sin(\pi \frac {t_i}q), \sin(\pi \frac {q-1-t_i}q\} } \bigg\}  .
\end{equation}
Now 
\begin{align*} 
\sum_{t=0}^{q-1} &  \min\bigg\{ q-1, 1+ \frac 1{  \min\{ \sin(\pi \frac {t}q), \sin(\pi \frac {q-1-t}q\} } \bigg\}\\
&= 3q-4 + 2 \sum_{1\leq t< \frac {q-1}2}  \frac { 1  }{  \sin(\pi \frac {t}q)} + \frac {  1_{2|q-1}  }{  \sin(\pi \frac {q-1}{2q})}
\end{align*}
The value of this sum is $\frac 2\pi q\log q+O(q)$, but for our application  we need the much weaker but fully explicit upper bound
\[
\leq (q-1) q^\tau \text{ for all } q\geq 133359
\]
where $\tau=\frac15-\eta$ and $\eta=10^{-9}$. The exponent ``$\frac 15$''  here is critical  because of the 
  $N^{4/5}$ in \eqref{eq: ExpSumPrime}). Substituting this into \eqref{eq: s-sum} we deduce that
\begin{equation} \label{eq: Full}
\sum_{j=0}^{q^k-1}  \bigg| S_{\mathcal A} \bigg(\theta +\frac j{q^k} \bigg)  \bigg|  \leq (q-1)^kq^{k\tau}.
\end{equation} 
Therefore the average value of $|S_{\mathcal A}(\alpha)|$ is given by
\begin{equation} \label{eq: FullIntegf}
\int_0^1 |S_{\mathcal A}(\alpha)|d\alpha = \int_0^{q^{-k}} \sum_{j=0}^{q^k-1}  \bigg| S_{\mathcal A} \bigg(\theta +\frac j{q^k} \bigg)  \bigg| d\theta
\leq \big( \tfrac{q-1}q \big)^k\cdot q^{k\tau}.
\end{equation}
This is  $<q^{k/5}=N^{1/5}$ much smaller than the   $N^{1/2}$ obtained from the mean square which is what is important in this argument. But it is also much  larger than $(C+o(1))^k$, the bound we obtained for $|S_{\mathcal A}(\theta)|$ for the typical   $\theta$ (that is, $\theta$ for which their base-$q$ digits are equidistributed) and it is feasible one can end up doing significantly better than we do here with cleverer arguments better exploiting the typical $\theta$.

\subsection{The mean value of $|S_{\mathcal A}'  (\alpha)|$} \label{sec 6.4}
For $n= \sum_{j=0}^{k-1} a_j q^j$ we have 
\[  
\frac d{d\theta} e(n\theta) = 2i\pi \cdot n e(n\theta) =2i\pi \cdot  \sum_{j=0}^{k-1} a_j q^j e(a_j q^j) \prod_{i\ne j} e(a_i q^i).
\]
We can modify the above argument from bounds for a sum of $ |S_{\mathcal A}(\cdot)|$-values to a sum of $ |S_{\mathcal A}'(\cdot)|$-values, by  bounding the contribution of the $j$th term in the product by $q^j$ times
\begin{align*}
 \bigg| \sum_{a=0}^{q-1} a\, e(a \phi)  - b\,  e(b \phi) \bigg|  &\leq \min\bigg\{ \frac{q(q-1)}2, b+ \bigg|  \frac{\sum_{j=1}^{q-1} e(j \phi)-(q-1)e(q \phi)}{ 1-e(\phi)}  \bigg| \bigg\} \\
 & \leq (q-1) \min\bigg\{ \frac{q}2, 1+   \frac{1}{ \sin(\pi \| \phi\| )}   \bigg\}\\
& \leq (q-1) \min\bigg\{ q-1, 1+   \frac{1}{ \sin(\pi \| \phi\| )}   \bigg\}
\end{align*}
 with $\phi=q^j\theta$.
Therefore, as $(q-1)\sum_{j=0}^{k-1}  q^j <q^k$ we obtain
\begin{equation} \label{eq: FullIntegf'}
\int_0^1 |S_{\mathcal A}'(\alpha)|d\alpha   \leq 2\pi (q-1)^k q^{k\tau}.
\end{equation}

\subsection{Bounds on $|S_{\mathcal A} (\theta_i)|$ at well spread-out points}

One can bound a differentiable function $f(\cdot)$ at a point by its values in a neighbourhood by the classical inequality
\[
|f(\theta) | \leq \frac 1{2\Delta} \int_{\theta-\Delta}^{\theta+\Delta} |f(\phi) |  d\phi + \frac 1{2} \int_{\theta-\Delta}^{\theta+\Delta} |f'(\phi) |  d\phi
\]
We can sum this over a set of points (on the unit circle), $\theta_1,\dots,\theta_m$ where  $|\theta_i-\theta_j|\geq 2\Delta$ if $i\ne j$ so the integrals above do not overlap, to obtain
\begin{equation} \label{spread}
\sum_{i=1}^m |f(\theta_i) |  \leq \frac 1{2\Delta} \int_0^1 |f(\phi) |  d\phi + \frac 1{2} \int_0^1 |f'(\phi) |  d\phi .
\end{equation}
Our choice of points is a bit complicated:  The $\theta_i$ will be selected within $\Delta= \frac 1{4S^2}$ of the fractions $\frac rs$ with $(r,s)=1$ and $0\leq r<s\leq S$ with $(r,s)=1$ displaced by a fixed quantity $\xi$.  The fractions are distinct so any two differ by 
$|\frac rs-\frac{r'}{s'}|\geq \frac{1}{ss'}> \frac{1}{S^2}$, and therefore the points differ by $\geq \frac 1{S^2}-2\Delta = 2\Delta$ and so
\[
\sum_{s\leq S} \sum_{\substack{0\leq r<s\\ (r,s)=1}} \max_{|\eta|\leq \Delta} \bigg|f\bigg(  \frac rs+\xi+\eta \bigg) \bigg|
\leq 2S^2 \int_0^1 |f(\phi) |  d\phi + \frac 1{2} \int_0^1 |f'(\phi) |  d\phi .
\]
We now apply this with $f=S_A$ and use \eqref{eq: FullIntegf} and \eqref{eq: FullIntegf'}   to obtain
\begin{equation} \label{eq: LS}
 \sum_{\substack{0\leq r<s\leq S \\ (r,s)=1}} \max_{|\eta|\leq  \frac 1{4S^2}} \bigg| S_A\bigg(  \frac rs+\xi+\eta \bigg) \bigg|
\leq (2S^2q^{-k} +\pi)(q-1)^k q^{k\tau}.
\end{equation}

  \subsection{Hybrid estimate}  We need notation that reflects that our sum is up to $q^k$, since we will now vary $k$. So let
  \[
  \widehat{A_k}(\theta) := S_{\mathcal A}(\theta) =\sum_{n\in \mathcal A(q^k)}  e(n\theta)
  \]
  Our formula \eqref{2}, implies that if $\ell\leq k$ then
  \[
\widehat{A_k}(\theta) = \widehat{A_{k-\ell}}(\theta) \widehat{A_\ell}(q^{k-\ell}\theta)  .
  \]
  For $m\leq k-\ell$ replace $k$ by $k-\ell$ and $k-\ell$ by $m$ so that
    \[
\widehat{A_{k-\ell}}(\theta) = \widehat{A_{m}}(\theta) \widehat{A_{k-\ell-m}}(q^{m}\theta)  ,
  \]
and therefore 
\[
\widehat{A_k}(\theta) =\widehat{A_{m}}(\theta) \widehat{A_{k-\ell-m}}(q^{m}\theta)   \widehat{A_\ell}(q^{k-\ell}\theta) .
\]
Since $|\widehat{A_{k-\ell-m}}(q^{m}\theta)|\leq (q-1)^{k-\ell-m}$ this yields
\[
|\widehat{A_k}(\theta)| =   (q-1)^{k-\ell-m} |\widehat{A_{m}}(\theta) |\cdot | \widehat{A_\ell}(q^{k-\ell}\theta) |,
\]
and so  
  \begin{align*}
\bigg| \widehat{A_k}  \bigg( \frac j {q^k} \bigg)\bigg|  &\leq (q-1)^{k-\ell-m}  \bigg| \widehat{A_{m}}\bigg( \frac j {q^{k}} \bigg) \bigg| \cdot \bigg| \widehat{A_\ell}\bigg( \frac j {q^{\ell}} \bigg) \bigg|\\
& \leq (q-1)^{k-\ell-m} \bigg| \widehat{A_\ell}\bigg( \frac j {q^{\ell}} \bigg) \bigg|  \cdot
\max_{i: |i-q^k \cdot \frac rs|\leq B }   \bigg| \widehat{A_{m}}\bigg( \frac i {q^{k}} \bigg) \bigg|  .
  \end{align*}
provided  $|j-q^k \cdot \frac rs|\leq B$.
 
 We let $B=q^\ell$ and $S^2=q^m$ so that $ 2S^2/q^m+\pi \ll 1$ and  $q^m=S^2\leq SM\ll N/B\ll q^{k-\ell}$.
 We have
\begin{align*}
\sum_{s\leq S} \sum_{\substack{0\leq r<s \\ (r,s)=1}} & \sum_{j: |j-q^k \cdot \frac rs|\leq B} \bigg| S_{\mathcal A}  \bigg( \frac j {q^k} \bigg)\bigg|      \\
\leq  (q-1)^{k-\ell-m} &   \sum_{\substack{0\leq r<s\leq S \\ (r,s)=1}}   \max_{i: |i-q^k \cdot \frac rs|\leq B }   \bigg| \widehat{A_{m}}\bigg( \frac i {q^k} \bigg) \bigg| \cdot \sum_{j: |j-q^k \cdot \frac rs|\leq B}  \bigg| \widehat{A_\ell}\bigg( \frac j {q^{\ell}} \bigg) \bigg|  .
  \end{align*}
  We extend the final sum to a sum over all $j \pmod {q^\ell}$ so that 
  \[
\sum_{j: |j-q^k \cdot \frac rs|\leq B}  \bigg| \widehat{A_\ell}\bigg( \frac j {q^{\ell}} \bigg) \bigg|  
\leq (q-1)^\ell q^{\tau \ell}
  \]  
by \eqref{eq: Full}, and therefore
\[
  \sum_{\substack{0\leq r<s\leq S \\ (r,s)=1}}  \sum_{j: |j-q^k \cdot \frac rs|\leq B} \bigg| S_{\mathcal A}  \bigg( \frac j {q^k} \bigg)\bigg|  \leq  (q-1)^{k-m}  q^{\tau \ell}  \sum_{\substack{0\leq r<s\leq S \\ (r,s)=1}}   \max_{i: |i-q^k \cdot \frac rs|\leq B }   \bigg| \widehat{A_{m}}\bigg( \frac i {q^k} \bigg) \bigg| .
\]
  
  For the next sum  we use that $B\leq  N/q^2SM$ and $S\leq M/q$ so that $B/N\leq 1/q^3S^2$. Therefore
\[
 \max_{i: |i-q^k \cdot \frac rs|\leq B }   \bigg| \widehat{A_{m}}\bigg( \frac i {q^k} \bigg) \bigg|  \leq 
 \max_{i: |\eta|\leq \frac{B}{q^k} }   \bigg| \widehat{A_{m}}\bigg( \frac rs + \eta \bigg) \bigg| \leq 
 \max_{i: |\eta|\leq \frac{1}{4S^2} }   \bigg| \widehat{A_{m}}\bigg( \frac rs + \eta \bigg) \bigg|
\]  
and so  the internal sum above is
\begin{align*}
\sum_{\substack{0\leq r<s\leq S \\ (r,s)=1}}   \max_{i: |i-q^k \cdot \frac rs|\leq B }   \bigg| \widehat{A_{m}}\bigg( \frac i {q^k} \bigg) \bigg| & \leq \sum_{\substack{0\leq r<s\leq S \\ (r,s)=1}}  \max_{i: |\eta|\leq \frac{1}{4S^2} }   \bigg| \widehat{A_{m}}\bigg( \frac rs + \eta \bigg) \bigg| \\ &\ll q^{O(1)} (q-1)^mq^{\tau m}
\end{align*}
by \eqref{eq: LS}. Therefore
\[
\sum_{s\leq S} \sum_{\substack{0\leq r<s \\ (r,s)=1}}  \sum_{j: |j-q^k \cdot \frac rs|\leq B} \bigg| S_{\mathcal A}  \bigg( \frac j {q^k} \bigg)\bigg|     
\ll_q |\mathcal A(N)|    q^{(\ell+m)\tau}.
\]
 which implies that \eqref{KeyMinors} holds as $q^{\ell+m}=BS^2$.

 \section{Reducing $q$}
 
 We have proved Maynard's Theorem, for primes missing one digit in base $q$, for all $q\geq 133359$.
 The goal is base $q=10$, so we need to find ways to improve the above argument to significantly reduce the size of $q$ to which it applies.       
 
 \subsection{More calculation}
 Now that we only have to work with the finite set of integers $q<133359$, and the finite set of values $b\in [0,q-1)$ we can do a separate calculation tailored more carefully to each individual case.
For example, instead of using  the bound $F_{\mathcal D}(\phi)\leq 1+ \frac 1{ \sin(\pi \| \phi \|)}$ we might instead work with the definition of $F_{\mathcal D}$ so that if $\phi\in [\frac tq,\frac{t+1}q)$with $t\in \mathbb Z$  then 
\[
 F_{\mathcal D}(\phi)\leq \max_{0\leq \eta<1} \bigg| \frac{e(\eta) -1} {e(\frac{t+\eta}q) -1}  - e(b\cdot \tfrac{t+\eta}q) \bigg| .
\]
Therefore we can replace the calculation after \eqref{eq: s-sum}, bounding the sum for each $i$, by the more precise
\[
\max_{0\leq b\leq q-1} \sum_{t=0}^{q-1} \max_{0\leq \eta<1} \bigg| \frac{e(\eta) -1} {e(\frac{t+\eta}q) -1}  - e(b\cdot \tfrac{t+\eta}q) \bigg|.
\]
For example if $q=101$, this improves the previous bound of $\leq 602.82\dots$ to something like $\leq 497$, but requires substantially more calculation.  Using this type of bound one gets weaker bounds for some $b$-values than for others, for a given $q$, and this ends up requiring more elaborate though stronger arguments.

 \subsection{A new cancelation}
 
 By \eqref{1} we have
 \[
 | S_{\mathcal A}(\theta)| \leq   \prod_{i=0}^{k-1} \min\bigg\{ q-1,   1+\bigg| \frac{e(q^{i+1}\theta) -1} {e(q^i\theta) -1}   \bigg|    \bigg\} .
  \] 
The second bound, $\leq 1+\frac 1{\sin( \pi \| q^i\theta\| ) }$, gives the minimum  if $1\leq t_i\leq q-2$.

 In  section \ref{sec: Minor} we proceeded by bounding the $i$th term of the product on average for each $i$, treating different $i$ independently (and so our upper bounds give  the ``worst case'' for each $i$).
We did so by simply using that $q^i\theta \mod 1\in [ \frac{t_i}q,\frac{t_i+1}q)$, and bounding 
$| e(q^{i+1}\theta) -1|\leq 2$.

This ignored the fact that 
$ \| q^{i+1}\theta \| $ can be determined given $ \| q^i\theta \| $. 
If we use the more precise 
 $q^i\theta \mod 1\in [ \frac{t_i+\frac{t_{i+1}}{q}}q, \frac{t_i+\frac{t_{i+1}+1}{q}}q)$ then 
   the upper bound \eqref{eq: sinbound} for the $i$th  and $(i+1)$st terms are
 \[
\leq  1+ \frac 1{ \sin(\pi \| \frac{t_i+\frac{t_{i+1}+\dots }{q}}q \|)} \text{ and }
 \leq     1+ \frac 1{ \sin(\pi \| \frac{t_{i+1}+\dots }{q} \|)}
 \]
 respectively, which are not independent but the dependence here is not so complicated, and we will be able to work with  this level of dependence.
 
 The idea is that we will obtain better upper bounds on $|S_{\mathcal A}(\theta)|$ by   taking each two consecutive terms of the product together. For example, 
\[
 |S_{\mathcal A}(\theta)| \leq  q   \prod_{j=0}^{k/2-1} R_{2j}
\]
where we take the $i$th and $(i+1)$st terms together, and
\[
R_i= \min\bigg\{ q-1,   1+\bigg| \frac{e(q^{i+1}\theta) -1} {e(q^i\theta) -1}   \bigg|    \bigg\} \cdot
\min\bigg\{ q-1,   1+\bigg| \frac{e(q^{i+2}\theta) -1} {e(q^{i+1}\theta) -1}   \bigg|    \bigg\}.
\]
Now if $1\leq t_i,t_{i+1}\leq q-2$ then 
 \begin{align*}
R_i&\leq \bigg( 1+\bigg| \frac{e(q^{i+1}\theta) -1} {e(q^i\theta) -1}   \bigg) \cdot
 \bigg( 1+\bigg| \frac{e(q^{i+2}\theta) -1} {e(q^{i+1}\theta) -1}   \bigg) \\
&\leq  1 + \bigg| \frac{e(q^{i+2}\theta) -1} {e(q^i\theta) -1}  \bigg|+ 
  \bigg| \frac{e(q^{i+1}\theta) -1} {e(q^i\theta) -1}   \bigg|  +\bigg| \frac{e(q^{i+2}\theta) -1} {e(q^{i+1}\theta) -1}   \bigg|  \\
  & \leq 1+  \frac{2+|e(q^{i+1}\theta) -1|} {|e(q^i\theta) -1|}    +  \frac{2} {|e(q^{i+1}\theta) -1|}   \\
  & \leq     1+ \frac { 1+ \max\{ \sin(\pi \frac {t_{i+1}+1}q)), \sin(\pi \frac {q-t_{i+1}}q)\} }{  \min\{ \sin(\pi \frac {t_i}q), \sin(\pi \frac {q-1-t_i}q)\} } 
+ \frac 1{  \min\{ \sin(\pi \frac {t_{i+1}}q), \sin(\pi \frac {q-1-t_{i+1}}q)\} }  .
\end{align*}
Summing our bounds over $0\leq t_i,t_{i+1}\leq q-1$ (using the upper bound $q-1$ on the $i$th term whenever $t_i$ equals $0$ or $q-1$, and similarly for the $(i+1)$st term) we get
\[
(3q-4)^2+ 
\bigg(  2 \sum_{1\leq t< \frac {q-1}2}  \frac { 1  }{  \sin(\pi \frac {t}q)} + \frac {  1_{2|q-1}  }{  \sin(\pi \frac {q-1}{2q})} \bigg) \bigg(   6q-8 +2 \sum_{2\leq u\leq  \frac {q}2}    \sin(\pi \tfrac {u}q) + \  1_{2|q-1}    \sin(\pi \tfrac {q+1}{2q})  \bigg) 
\]
which is $<(q-1)^2q^{2/5}$ for $q\geq 18647$ (by a computer calculation), and therefore we have proved the claimed result for such $q$.

We can combine this improvement with that of the previous subsection and the two ideas together should improve the bound on $q$ further.

By taking two consecutive $i$-values together we have improved our lower bound on $q$ by   factor of more than 7, so we can probably get further improvements if  we multiply together three consecutive $i$-values, or more? When we  do this, it is natural to ask how  to keep track of useful cancelations, like the
\[
 \frac{e(q^{i+2}\theta) -1} {e(q^{i+1}\theta) -1}  \cdot  \frac{e(q^{i+1}\theta) -1} {e(q^i\theta) -1}  =  \frac{e(q^{i+2}\theta) -1} {e(q^i\theta) -1} 
\]
used above, and when do we chose to use the trivial upper bound ``2'' on the numerator?  Maynard's surprising idea is to keep track of all this by regarding the different terms of the product, averaged over all possible sets of $t_i$'s, as transition probabilities in a Markov process.

\subsection{Better bounds on \eqref{2} via a special Markov process} 
We approximated the terms in \eqref{2} using only the first term of the base-$q$  expansion of $q^i\theta \mod 1$. However if we obtain a more precise approximation using, say, the first two terms, $t_i$ and $t_{i+1}$, of the base-$q$  expansion of $q^i\theta \mod 1$, then the bounds on the $i$th and $(i+1)$st terms are no longer independent (it was that independence which allowed us  to take a  sum of the product   equal to the product of various smaller sums).  In particular we obtain a more accurate approximation using 
  $e(q^i\theta)\approx e(t_i/q+ t_{i+1}/q^2)$ which involves the first two terms of the expansion.  Substituting this  approximation into \eqref{2} 
 yields that   
 \[
 | S_{\mathcal A}(\theta)| \approx \prod_{i=0}^{k-1}  F(t_i,t_{i+1}) \text{ where } F(t,u):=\bigg|\frac{e(\frac uq) -1} {e(\frac tq +\frac u{q^2} ) -1}  - e(b (\tfrac tq +\tfrac u{q^2} ) )\bigg| \text{ if } t\ne 0,
 \]
 and $F(0,u)=q-1$.   Now the consecutive terms depend on each other so we cannot separate them as before. Instead we can
 form the $q$-by-$q$ matrix $M$ with entries $M_{a,b}:=\frac{F(a,b)}{q-1}$ for $0\leq a,b\leq q-1$.  Then for $t_0,t_k\in \{ 0,\dots,q-1\}$
 \[
 (q-1)^kM^k_{t_0,t_k} = \sum_{t_1,\dots ,t_{k-1}\in \{ 0,\dots,q-1\} }     \prod_{i=0}^{k-1}  F(t_i,t_{i+1})
 \approx  \sum_{\substack{t_1,\dots ,t_{k-1} \in \{ 0,\dots,q-1\} \\ \theta=\sum_{i=0}^{k} t_i/q^{i+1}  }}  | S_{\mathcal A}(\theta)|.
 \] 
  Summing this over all $t_0,t_k\in \{ 0,\dots,q-1\}$ gives the complete sum over the $\theta=j/q^k$; that is,
  \[
 (q-1)^{-k} \sum_{j=0}^{q^k-1}  | S_{\mathcal A}(\tfrac j{q^k})| \approx (1,1,\dots, 1) M^k (1,1,\dots, 1)^T \leq c_M |\lambda_M|^k
  \]
  where $\lambda_M$ is the   largest eigenvalue of $M$ and $c_M>0$ is some computable constant.\footnote{We need to change the ``$\approx$'' in $| S_{\mathcal A}(\theta)| \approx \prod_{i=0}^{k-1}  F(t_i,t_{i+1}) $ above to a precise inequality, like
  \[
| S_{\mathcal A}(\theta)| \leq \prod_{i=0}^{k-1}  F(t_i,t_{i+1}),  \text{ where }   F(t,u):=\max_{0\leq \eta \leq 1/q^2} \bigg|\frac{e(\frac {u+\eta}q) -1} {e(\frac tq +\frac {u+\eta}{q^2} ) -1}  - e(b (\tfrac tq +\tfrac {u+\eta}{q^2}  ) )
  \bigg| 
  \] }
  Our  proof of the bounds for the minor arcs can be   modified in a straightforward way, and then the result follows provided 
  \[
  \lambda_M< q^{1/5}.
  \]
 With our earlier proved results we can assume that  $q< 18647$; in particular we can compute the matrix in each case and determine the largest eigenvalue.
 
\subsection{A more general Markov process}
  But this is far from the end of the story, since we can be more precise by replacing the transition from the first two terms of the expansion of $q^i\theta$,
  $t_i,t_{i+1}$, to the next two, $t_{i+1},t_{i+2}$, in our ``Markov process'', by the transition from the first $\ell$ terms of the expansion of $q^i\theta$ to the next $\ell$. This yields a $q^\ell$-by-$q^\ell$ transition matrix $M=M^{(\ell)}$ which is indexed by $\ell$ digits in base $q$ and   $(M^{(\ell)})_{I,J}$ can only be non-zero if
\[
I=(t_1,\dots,t_\ell), J = (t_2,\dots, t_{\ell+1}) \text{ for some base-$q$ digits } t_1,\dots,t_{\ell+1}.
\]
Therefore each row and column is supported at only $q$ entries. 
 
 If $\theta=\sum_{i=1}^{\ell+1} t_i/q^i$ then the corresponding entry of $M^{(\ell )}$ is $G(t_1,\dots, t_{\ell+1}) $ where 
\[
G(t_1,\dots, t_{\ell+1}):= \max_{0\leq \eta\leq 1/q^{\ell+1}}\frac 1{|\mathcal D|} F_{\mathcal D}(\theta+\eta) .
\] 
If $\lambda_{\ell}$ is the largest eigenvalue of $M^{(\ell)}$ in absolute value then 
\[
  \sum_{j=0}^{q^k-1}  \bigg| S_{\mathcal A}\bigg(\frac j{q^k}\bigg)\bigg|  \ \ll |\mathcal D|^k\cdot  |\lambda_{\ell }|^{k} ,
\]
and therefore if $|\lambda_{\ell }|< q^{1/5}$ for some $\ell\geq 1$ then there are indeed the expected number of primes with base-$q$ digits in the set $\mathcal D$. 

Since these are truncations of the true Markov process on a Hilbert space (with infinitely many dimensions) we have that  $|\lambda_1|>|\lambda_2|>\dots$ and so our bounds improve as $\ell$ gets larger. These tend to a  (positive) limit $|\lambda_\infty|$
which gives the solution to the eigenvalue problem for the matrices in this Hilbert space.
However, numerical approximation shows that $|\lambda_\infty|$ is not as small as would be needed to resolve the base 10 problem.

\subsection{Using the Markov process to remove generic minor arcs in small bases}  

Maynard's next idea for small $q$ was to ``remove'' as many ``generic'' minor arcs as possible. He does so by using a simple moment argument: For any $\sigma>0$ we have 
\begin{equation} \label{eq: moments}
\#\bigg\{ j\in [0,N): \bigg| S_{\mathcal A}\bigg(\frac j{q^k}\bigg)\bigg| \geq \frac{A(N)}T\bigg\} \leq
\bigg(  \frac T{A(N)} \bigg) ^\sigma  \sum_{j=0}^{N-1}  \bigg| S_{\mathcal A}\bigg(\frac j{q^k}\bigg)\bigg|^\sigma ,
\end{equation}
so now we are interested in bounding the $\sigma$th moment of $|S_A|$. To do this we work with the matrix 
$M^{(\ell,\sigma)}$ where $ (M^{(\ell,\sigma)})_{I,J}=(M^{(\ell)})_{I,J}^\sigma$, so that if $\lambda_{\ell,\sigma}$ is the largest eigenvalue of $M^{(\ell,\sigma)}$ in absolute value then 
\[
\frac 1{A(N)^\sigma}  \sum_{j=0}^{q^k-1}  \bigg| S_{\mathcal A}\bigg(\frac j{q^k}\bigg)\bigg|^\sigma \ \ll |\lambda_{\ell,\sigma}|^{k} .
\]
On the other hand
\begin{equation*}
\#\bigg\{ j\in [0,N): \bigg| S_{\mathcal P}\bigg(\frac j{q^k}\bigg)\bigg| \geq U\bigg\} \leq
U^{-2}  \sum_{j=0}^{N-1}  \bigg| S_{\mathcal P}\bigg(\frac j{q^k}\bigg)\bigg|^2 = U^{-2} N \pi(N)\sim \frac{N^2}{U^2\log N}.
\end{equation*}
Therefore if $\mathcal E=\{ j\in [0,N): |S_{\mathcal A} (\frac j{q^k}) | \geq \frac{A(N)}T \text{ or } 
 |S_{\mathcal P} (\frac j{q^k}) | \geq U\}$ then
 \[
 \frac 1N\sum_{\substack{j=0 \\  j\not\in \mathcal E}}^{N-1} \bigg| S_{\mathcal A}\bigg(\frac j{q^k}\bigg) \cdot S_{\mathcal P}\bigg(\frac j{q^k}\bigg)\bigg|  \leq \frac{A(N)}{(\log N)^2}
 \]
taking $U=T/(\log N)^2$. Now if $|\lambda_{\ell,\sigma}|< q^\rho$ then 
\[
|\mathcal E|\ll T^\sigma N^\rho + \frac{N^2(\log N)^3}{T^2}< N^{\frac{2+\rho+\rho\sigma}{2+\sigma}+o(1)},
\]
selecting $T=N^{\frac{2-\rho}{2+\sigma}}$.

Karwatowski \cite{Kar} used the fact that the eigenvalues of a matrix are bounded in absolute value by the largest sum of the absolute values of the elements in a row of the matrix, to numerically prove the bounds
 \[
\lambda_{4,1}< q^{\frac {27}{77}} \text{ and } \lambda_{4,\frac{235}{154}}< q^{\frac {59}{433}}
\]
for all $q\geq 10$ 
(Maynard had already shown   these inequalities hold for $q=10$.) 
The moment method with $\sigma=\frac{235}{154}$ then implies $|\mathcal E|\ll N^{2/3}$ arguing as above, and therefore one can focus on the exceptional $j$-values.

To make the base-$10$ argument unconditionally doable, Maynard developed delicate sieve methods. In effect this allowed him to replace needing to understand how often primes are written with the digits from $\mathcal D$ in base $q$, to understanding when integers  composed of a product of a few large  primes in certain given intervals are so represented.   Maynard could therefore   improve on the upper bounds for exponential sums over primes (as in \eqref{eq: ExpSumPrime}) when appropriately weighted, since now he was working with a more malleable set of the integers. He was able to restrict attention to a set $\mathcal E\subset \mathfrak m$ of exceptional integers $j$ with $|\mathcal E|\ll N^{.36}$.

\subsection{The exceptional minor arcs}  
If  $j/q^k\in \mathcal E$  has an important effect on our sum, then  the fraction $j/q^k$ will have to   simultaneously have several surprising Diophantine features, which Maynard proves are mostly incompatible (when $q=10$). The techniques are too complicated to discuss here. The following diagram exhibits the tools used in the whole proof, but especially when dealing with these exceptional arcs.

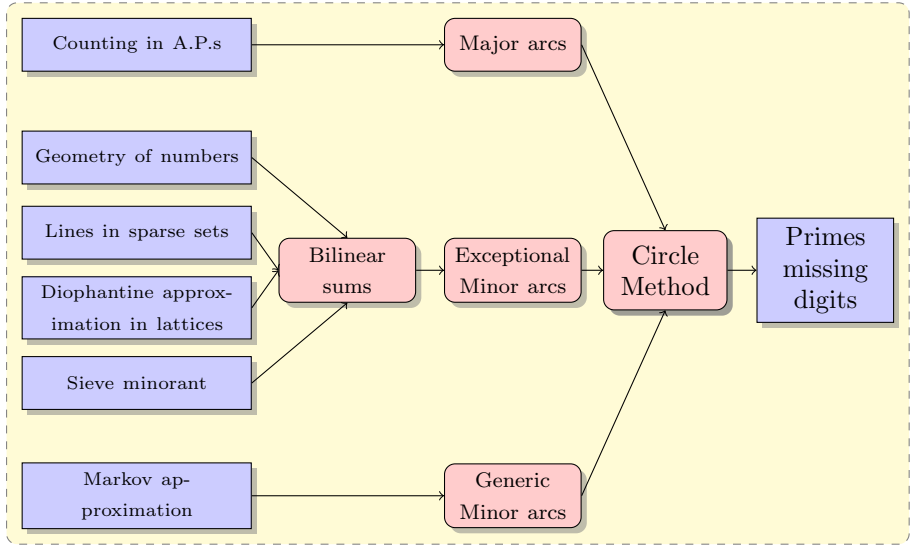
\begin{figure}[H]
   \begin{center}
\makebox[\textwidth]{\parbox{1.5\textwidth}{
\begin{center}
   \pgfdeclarelayer{background}
   \pgfdeclarelayer{foreground}
   \pgfsetlayers{background,main,foreground}
 
   \tikzstyle{interface}=[draw, fill=blue!20, text width=8em,
      text centered, minimum height=2.0em,drop shadow]
   \tikzstyle{daemon}=[draw, text width=4.5em, fill=red!20,
      minimum height=2em, text centered, rounded corners, drop shadow]
   \tikzstyle{dots} = [above, maximum width=7em, text centered]
   \tikzstyle{wa} = [daemon, text width=4em, fill=red!20,
      minimum height=3em, rounded corners, drop shadow]
   \tikzstyle{ur}=[draw, fill=blue!20, text width=4.5em,
      text centered, minimum height=2.0em,drop shadow]
 
   \def\blockdist{1.3}
   \def\edgedist{0.}

   \begin{tikzpicture}
 
      \node (wa) [wa]  {Circle\\ Method};
      \path (wa.west)+(-6.2,3) node (i1)[interface] {\scriptsize Counting in A.P.s};
      \path (wa.west)+(-6.2,1.5) node (i2)[interface] {\scriptsize Geometry of numbers};
      \path (wa.west)+(-6.2,0.5) node (i3)[interface] {\scriptsize Lines in sparse sets};
      \path (wa.west)+(-6.2,-0.5) node (i6)[interface] {\scriptsize Diophantine approximation in lattices};
      \path (wa.west)+(-6.2,-1.5) node (i5)[interface] {\scriptsize Sieve minorant};
      \path (wa.west)+(-6.2,-3) node (i4)[interface] {\scriptsize Markov approximation};
      \path (wa.west)+(-3.4,0.0) node (d5)[daemon] {\footnotesize Bilinear\\ sums};
      \path (wa.west)+(-1.2,+3) node (d1)[daemon] {\footnotesize Major arcs};
      \path (wa.west)+(-1.2,-3) node (d2)[daemon] {\footnotesize Generic Minor arcs};
      \path (wa.east)+(1.3,0.0) node (ur)[ur] {Primes missing digits};
      \path (wa.west)+(-1.2,0.0) node (d3)[daemon] {\footnotesize Exceptional Minor arcs};
      \path [draw, ->] (wa.east) -- node [above] {} (ur.west) ;
      \path [draw, ->] (d1.east) -- node [above] {} (wa.north) ;
      \path [draw, ->] (d2.east) -- node [above] {} (wa.south) ;
       \path [draw, ->] (d3.east) -- node [above] {} (wa.west) ;
       \path [draw, ->] (i1.east) -- node [above] {} (d1.west) ;
        \path [draw, ->] (i4.east) -- node [above] {} (d2.west) ;
        \path [draw, ->] (i5.east) -- node [above] {} (d5.south) ;
                \path [draw, ->] (i6.east) -- node [above] {} (d5.west) ;
          \path [draw, ->] (i2.east) -- node [above] {} (d5.north) ;
           \path [draw, ->] (i3.east) -- node [above] {} (d5.west) ;
                                  \path [draw, ->] (d5.east) -- node [above] {} (d3.west) ;
 
      \begin{pgfonlayer}{background}
            \path (i2.west |- i1.north)+(-0.2,0.2) node (a) {};
            \path (ur.east |- i4.south)+(+0.2,-0.2) node (c) {};
 
            \path[fill=yellow!20,rounded corners, draw=black!50, dashed]
                  (a) rectangle (c);
      \end{pgfonlayer}
 
   \end{tikzpicture}
\end{center}}}
   \end{center}
   \vspace{-0.5cm}
\caption{Outline of steps to prove primes with missing digits}
\end{figure}

\section{Generalizations}  
Our argument for sufficiently large $q$, generalizes to a given set $\mathcal D$, if 
$\mathcal D$ contains two consecutive integers (for section \ref{sec 5.3}), and  if
\[
\sum_{t=0}^{q-1}  \max_{0\leq \eta < \frac 1q} F_{\mathcal D}\bigg(\frac tq +\eta \bigg) < (q-|\mathcal D|) q^{1/5}
\]
(The contributions of the $\sum_{a\in \mathcal D} ae(a\phi)= q^{O(1)}$ in  section \ref{sec 6.4} and so are not relevant).  Now if $\mathcal D = \{ 0,\dots ,q-1\} \setminus \mathcal R$ for a set $\mathcal R$ with $r$ elements then 
\[
\bigg| \sum_{a\in \mathcal D}  e(a\phi)\bigg| \leq \bigg| \sum_{b\in \mathcal R}  e(b\phi)\bigg|  +  \bigg| \sum_{a=0}^{q-1}  e(a\phi)\bigg|   \leq r + \frac 1{\sin( \pi \| \phi\|)},
\]
so that 
\[
\sum_{t=0}^{q-1}  \max_{0\leq \eta < \frac 1q} F_{\mathcal D}\bigg(\frac tq +\eta \bigg)
\leq (q-1)r+q\log q+O(q) < (q-r) q^{1/5}
\]
if $r<(1-\epsilon) q^{1/5}$ for $q$ sufficiently large. 
We can improve this using \eqref{spread}, first summing over the points with $t$ even, then those with $t$ odd, this is
\begin{align*} 
  \leq q \int_0^1 \bigg| \sum_{b\in \mathcal R} e(b\phi) \bigg|  d\phi 
+  \int_0^1\bigg| \sum_{b\in \mathcal R}b e(b\phi) \bigg|  d\phi < 2q\sqrt{r}
\end{align*}
since, for any  coefficients $c_b$
\[
\bigg(  \int_0^1 \bigg| \sum_{b\in \mathcal R} c_b e(b\phi) \bigg|  d\phi \bigg)^2 \leq
\int_0^1 \bigg| \sum_{b\in \mathcal R} c_b e(b\phi) \bigg|^2  d\phi = \sum_{b\in \mathcal R} |c_b|^2
\]
by the Cauchy-Schwarz inequality. Therefore 

\emph{There are roughly the expected number of primes whose base-$q$ digits come from the set $\mathcal D$ whenever $|\mathcal D|\geq q-\frac{1}5q^{2/5}$, for $q$ sufficiently large.}
\
\bigskip

Another idea is to let $\mathcal D$ be a set of $r$ consecutive integers; we can see that
\[
\bigg| \sum_{a\in \mathcal D}  e(a\phi)\bigg| \leq \min \bigg\{ r, \frac 1{\sin( \pi \| \phi\|)} \bigg\}
\]
so that 
\[
\sum_{t=0}^{q-1}  \max_{0\leq \eta < \frac 1q} F_{\mathcal D}\bigg(\frac tq +\eta \bigg)
\ll (q-r) \frac qr \log r    
\]
and this is $<q^{1/5}$ provided $r\gg q^{4/5}\log q$. Therefore

\emph{There are roughly the expected number of primes whose base-$q$ digits come from any set $\mathcal D$ of  $\gg q^{4/5}\log q$ consecutive integers, for $q$ sufficiently large.}\smallskip

\noindent  The ``$\frac 45$'' was improved to ``$\frac 34$''   in  \cite{May5}, and even 
to  ``$\frac {57}{80}$'' if one just wants a lower bound of the correct order of magnitude.

\newpage

\part{Approximations by reduced fractions}
\section{Approximating most real numbers}
We begin by reducing the real numbers modulo the integers; that is, given $\theta\in \mathbb R$ we consider the equivalence class $(\theta)$ of real numbers that differ from $\theta$ by an integer (and so each such $(\theta)$ is represented by a unique  real number  in $(-\frac 12, \frac 12]$).

Dirichlet observed that if $\alpha\in [0,1)$ then the representations of 
\[
(0), (\alpha), (2\alpha),\cdots, (N\alpha)
\]
 all belong to an interval of length $1$ so two of them $(i \alpha)$ and $(j \alpha)$ must differ by   $<\frac 1N$, by the pigeonhole principle.\footnote{Moreoverm by embedding the interval onto the circle by the map $t\to e(t):=e^{2i\pi t}$ we see that they must differ by  $<\frac 1{N+1}$.}
Now if $n=|j-i|$ then $n\leq N$ and $n=\pm(j-i)$, so that
\[
\pm (n\alpha) \equiv \pm n\alpha = (j-i)\alpha \equiv (j\alpha) -  (i\alpha) \mod 1.
\]
Therefore there exists an integer $m$ for which $|n\alpha - m| < \frac 1N$ which we rewrite as
\[
\bigg| \alpha -\frac mn \bigg| < \frac 1{nN} \leq \frac 1{n^2}.
\]
This is a close approximation to $\alpha$ by rationals, and one wonders whether one can do much better?
In general, no, since the continued fraction of the golden ratio $\phi:=\frac{1+\sqrt{5}}2$ implies that the best approximations to $\phi$ are given by $F_{n+1}/F_n, n\geq 1$ where $F_n$ is the $n$th Fibonacci number: One can show that 
\[
\bigg| \phi -\frac {F_{n+1}}{F_n} \bigg| \sim \frac 1{\sqrt{5}} \cdot \frac 1{F_n^2},
\]
and so all approximations to $\phi$ by rationals $p/q$ miss by $\geq \{ 1+o(1)\} \frac 1{\sqrt{5}} \cdot \frac 1{q^2}$.

This led researchers at the end of the 19th century to realize that if the partial quotients in the continued fraction for irrational  $\alpha$ are bounded, say by $B$ (note that $\phi=[1,1,1,\dots]$) then there exists a constant $c=c_B>0$ such that 
$|\alpha-\frac mn|\geq \frac{c_B}{n^2}$. However there are very few such $\alpha$   under any reasonable measure. If the partial quotients aren't bounded then how good can approximations be? And how well can one approximate famous irrationals like $\pi$? (still a very open question).\footnote{If $\alpha$ has continued fraction $[a_0,a_1,\dots]$ and  
$|\alpha-\frac bq| <\frac 1{2q^2}$ then
$\frac bq$   a convergent of the continued fraction, say the $j$th convergent, and then  
$|\alpha-\frac {b_j}{q_j}| \asymp \frac 1{a_jq_j^2}$; that is, we get better approximations the larger the $a_j$ in the continued fractions (especially in comparison to the $q_j$).  
However we do not understand the continued fractions of most real numbers $\alpha$ well enough to be able to assert that the problem is resolved,  so we have  transferred the difficulty of the problem into a seemingly different domain. See appendix 11B of \cite{GrM} for more on continued fractions.}

An easy argument shows that the set of $\alpha\in [0,1)$ with infinitely many rational approximations $\frac mn$ for which
$|\alpha-\frac mn|\leq \frac 1{n^3}$ has measure $0$. Indeed if there are infinitely many such rational approximations then there is one with $n>x$ (an integer). Now for each $n$ the measure of $\alpha\in [0,1)$ with $|\alpha-\frac mn|\leq \frac 1{n^3}$ is $\frac 1{n^3}$ for $m=0$ or $n$, 
$\frac 2{n^3}$ for $1\leq m\leq n-1$ and $0$ otherwise, a total of $\frac 2{n^2}$, and summing that over all $n>x$ gives
$\sum_{n>x} \frac 2{n^2}<\int_x^\infty \frac 2{t^2} dt = \frac 2x$. Letting $x\to \infty$ we see that the measure is $0$.
Obviously the analogous result holds for $|\alpha-\frac mn|\leq \frac 1{(n\log n)^2}$, and any other such bounds that lead to convergence of the infinite sum.

More generally we should study, for a given function $\psi: \mathbb Z_{\geq 1}\to \mathbb R_{\geq 0}$,
 the set $\mathcal L(\psi)$ which contains those $\alpha\in [0,1)$ 
for which there are infinitely many rationals $m/n$ for which
\[
\bigg| \alpha -\frac mn \bigg| \leq \frac{\psi(n)}{n^2}.
\]
We have seen that $\mathcal L(1)=[0,1)$ whereas if $c< 1/\sqrt{5}$ then $\phi-1\not\in \mathcal L(c)$ so $\mathcal L(c)\neq [0,1)$.
Moreover if $\sum_n \psi(n)/n$ is convergent then $\mu(\mathcal L(\psi))=0$ where $\mu(\cdot)$ is the Lebesgue measure.
In each case that  we have worked out, $\mu(\mathcal L(\psi))=0$ or $1$, and Cassels \cite{Cas} showed that this is always true (using the Birkhoff Ergodic Theorem)!
So we need only decide between these two cases.

The first great theorem in \emph{metric Diophantine approximation} was due to Khinchin who showed that if $\psi(n)$ is a decreasing function then
\[
\mu(\mathcal L(\psi))= \begin{cases} 0\\ 1 \end{cases} \text{ if and only if } \sum_{n\geq 1} \frac{ \psi(n)}n \text{ is }
 \begin{cases}  \text{ convergent}\\  \text{ divergent} \end{cases}.
 \]
 Thus  measure $1$ of reals $\alpha$ have approximations $\frac mn$ with $| \alpha -\frac mn| \leq \frac 1{n^2\log n}$, and 
 measure $0$ with $| \alpha -\frac mn| \leq \frac 1{n^2(\log n)^{1+\epsilon}}$

The hypothesis ``$\psi(n)$ is  decreasing'' is too restrictive since, for example, one can't determine anything from this about rational approximations where the denominator is prime. So can we do without it? Our proof above that 
if $\sum_{n\geq 1} \frac{ \psi(n)}n $ is convergent then $\mu(\mathcal L(\psi))=0$,  works for general $\psi$.
Indeed we follow the usual proof of the first Borel-Cantelli lemma:
Let $E_n$ be the event that $\alpha\in [\frac mn -\frac{\psi(n)}{n^2},\frac mn +\frac{\psi(n)}{n^2}]\cap [0,1]$ for some $m\in \{ 0,1,\dots,n\}$, where we have selected $\alpha$ randomly from $[0,1]$, and we established that
$\sum_n \mathbb P(E_n)=\sum_{n} \frac{ \psi(n)}n<\infty$. Then, almost surely, only finitely many of the $E_j$ occur, and so $\mu(\mathcal L(\psi))= 0$.  

The second Borel-Cantelli lemma states that if the $E_n$ are independent and $\sum_n\mathbb P(E_n)$ diverges then almost surely  infinitely many of the $E_j$ occur.  Our $E_n$ are far from independent (indeed compare $E_n$ with $E_{2n}$) but this  nonetheless suggests that perhaps with the right notion of independence it is feasible that Khinchin's theorem holds without the decreasing condition. 

\subsection{Duffin and Schaefer's example}    Duffin and Schaefer constructed a (complicated) example of $\psi$ for which $ \sum_{n\geq 1} \frac{ \psi(n)}n $ diverges but $\mu(\mathcal L(\psi))= 0$;. Their example uses many representations like $\frac 13=\frac 26$, that is, non-reduced fractions:

We begin with $\psi_0$ where $\psi_0(q)=0$ unless $q=q_\ell:=\prod_{p\leq \ell} p$ is the product of the primes up to some prime $\ell$, in which case  $\psi_0(q_\ell)=\frac {q_\ell}{\ell \log \ell}$. Therefore 
\[
\sum_q  \frac{\psi_0(q)}q = \sum_\ell \frac {1}{\ell \log \ell} 
\]
which converges by the prime number theorem, and so $\mu(\mathcal L(\psi_0))= 0$ as we just proved in the last subsection.

Now we construct a new $\psi$ for which if $q$ is  squarefree integer  with largest prime factor $\ell$ (so that $q$ divides $q_\ell$), then $\psi(q)=q^2/(q_\ell \ell \log \ell)$,  and $\psi(q)=0$ otherwise. Now if  $|x-\frac aq|\leq \frac{\psi(q)}{q^2}$ then for $A=a(q_\ell/q)$ we have
\[
\bigg|x-\frac A{q_\ell}\bigg| =\bigg|x-\frac aq\bigg| \leq \frac{\psi(q)}{q^2} =  \frac{\psi(q_\ell)}{q_\ell^2}
=  \frac{\psi_0(q_\ell)}{q_\ell^2}
\] 
so that $\mathcal L(\psi)=\mathcal L(\psi_0)$ which has measure $0$. On the other hand
\[
\sum_q  \frac{\psi(q)}q = \sum_{\ell} \frac 1{\ell\log \ell} \sum_{\ell|q|q_\ell} \frac q{q_\ell}
=  \sum_{\ell} \frac 1{\ell\log \ell} \prod_{p<\ell}\bigg( 1+ \frac 1p\bigg) \gg  \sum_{\ell} \frac 1\ell
\]
by Mertens' Theorem, which diverges.

\subsection{A revised conjecture} 
Duffin and Schaefer's example  uses many representations like $\frac 13=\frac 26$, which suggests that we should restrict attention to \emph{reduced fractions} $\frac mn$ with $(m,n)=1$.
We let $E_n^*$ be the event that $\alpha\in [\frac mn -\frac{\psi(n)}{n^2},\frac mn +\frac{\psi(n)}{n^2}]\cap [0,1]$ for some $m\in \{ 0,1,\dots,n\}$ with $(m,n)=1$.

Therefore Duffin and Schaefer defined $\mathcal L^*(\psi)$ to be those $\alpha\in [0,1)$ 
with infinitely many reduced fractions $m/n$ for which
\[
\bigg| \alpha -\frac mn \bigg| \leq \frac{\psi(n)}{n^2},
\]
and conjectured 
\[
\mu(\mathcal L^*(\psi))= \begin{cases} 0\\ 1 \end{cases} \text{ if and only if } 
\sum_{n\geq 1}  \frac{ \phi(n) }n \cdot  \frac{ \psi(n) }n \text{ is }
 \begin{cases}  \text{ convergent}\\  \text{ divergent} \end{cases}.
 \]
 Here $\phi(n)=\# \{ \frac mn\in [0,1): (m,n)=1\}$. Now if
$\sum_n \mathbb P(E_n^*)=\sum_{n} \frac{ \phi(n) }n \cdot \frac{ \psi(n)}n<\infty$, then 
almost surely, only finitely many of the $E_j^*$ occur, and so $\mu(\mathcal L^*(\psi))= 0$.
We therefore can assume that 
$\sum_{n\geq 1}  \frac{ \phi(n) }n \cdot  \frac{ \psi(n) }n$ is divergent.

Gallagher  \cite{Gal} (in a slight variant of Cassell's result \cite{Cas}) showed that  $\mu(\mathcal L^*(\psi))$ always equals either $0$ or $1$. Therefore we only need to show that $\mu(\mathcal L^*(\psi))>0$ to deduce that   $\mu(\mathcal L^*(\psi))=1$.

 Duffin and Schaefer themselves proved the conjecture in the case that 
there are arbitrarily large $Q$ for which
\[
\sum_{q\leq Q}    \frac{ \phi(q) }q \cdot  \frac{ \psi(q) }q \gg \sum_{q\leq Q}    \frac{ \psi(q) }q;
\]
which more-or-less implies that the main weight of $\psi(q)$ should not be focussed on integers $q$ with many small prime factors (which are extremely rare), since that is what forces 
\[
\frac{\phi(q)}q = \prod_{p|q} \bigg( 1-\frac 1p\bigg) \text{ to be small.}
\] Thus for example, the conjecture follows if we only allow prime $q$ (that is, if $\psi(q)=0$ whenever $q$ is composite), or if we only allow integers $q$ which have no prime factors $<\log q$.

 In 2021, Koukoulopoulos and Maynard \cite{KM} showed that this Duffin-Schaefer conjecture is true, the end of a long saga.
 The proof is a blend of number theory, probability theory, combinatorics, ergodic theory, and graph theory combined with considerable ingenuity.

 \subsection{Probability} Assuming that
$\sum_{n\geq 1}  \frac{ \phi(n) }n \cdot  \frac{ \psi(n) }n$ is divergent, we want to show that 
almost surely, infinitely many of the $E_j^*$ occur, where $E_q^*$ is the event that 
$\alpha$ belongs to
 \[
  [0,1)\cap \bigcup_{(a,q)=1}  \bigg[  \frac aq - \frac{\psi(q)}{q^2}, \frac aq + \frac{\psi(q)}{q^2}\bigg].
 \]
 The $E_q^*$ are not ``independent'', but were they independent enough, say if
\[
\mu(E^*_q\cap E^*_r) = (1+o_{q,r\to \infty} (1)) \, \mu(E^*_q)  \, \mu(E^*_r) ,
\]
then we could prove our result; however  one can easily find counterexamples to this, for example when $r=2q$. 
On the other hand, since  we only need to show that $\mu(\mathcal L^*(\psi))>0$, we will only need to establish a very weak quasi-independence, on average, like 
\begin{equation}\label{eq: Assume 1}
  \sum_{Q\leq q\ne r<R} \mu(E^*_q\cap E^*_r)\leq  10^6  \bigg(\sum_{Q\leq q <R} \mu(E^*_q) \bigg)^2
\end{equation}
for arbitrarily large $Q$ and certain $R$: To prove this note that since $ \sum_{q\geq Q} \mu( E^*_q)= 2   \sum_{q\geq Q} \frac{ \phi(q) }q \cdot  \frac{ \psi(q) }q $   diverges,  we may select $R\geq Q$ for which $1\leq  \sum_{Q\leq q<R} \mu( E^*_q)\leq 2$.  
Now let $N=\sum_{Q\leq q<R} 1_{E^*_q}$ so that $\mathbb E[N]=\sum_{Q\leq q<R} \mu(E^*_q)$ and so
\begin{align*}  
1\leq \bigg( \sum_{Q\leq q<R} \mu(E^*_q) \bigg)^2&=   \mathbb E[N]^2= \mathbb E[1_{N>0}\cdot N]^2\leq 
  \mu\bigg( \bigcup_{Q\leq q<R} E^*_q\bigg)\cdot \mathbb E[N^2] \\
  &=  \mu\bigg( \bigcup_{Q\leq q<R} E^*_q\bigg) \sum_{Q\leq q,r<R}  \mu(E^*_q\cap E^*_r)  
\end{align*}  
by the Cauchy-Schwarz inequality. Therefore 
\[
\mu\bigg( \bigcup_{q\geq Q} E^*_q\bigg) \geq \mu\bigg( \bigcup_{Q\leq q<R} E^*_q\bigg) \geq 10^{-6}
\]
by \eqref{eq: Assume 1}. But this is true for arbitrarily large $Q$ and so  $\mu(\mathcal L^*(\psi))\geq 10^{-6}$, which implies that  $\mu(\mathcal L^*(\psi))=1$. 

 Following  Pollington and Vaughan \cite{PV} we study $\mu(E^*_q\cap E^*_r)$, assuming $(q,r)=1$ for convenience: If $\alpha\in  [  \frac aq - \frac{\psi(q)}{q^2}, \frac aq + \frac{\psi(q)}{q^2}]\cap
  [  \frac br - \frac{\psi(r)}{r^2}, \frac br + \frac{\psi(r)}{r^2}]$ with $(a,q)=(b,r)=1$ then
  $|\frac aq - \frac br|\leq \frac{\psi(q)}{q^2}+  \frac{\psi(r)}{r^2}\leq 2\Delta$ where $\Delta:=\max\{ \frac{\psi(q)}{q^2},\frac{\psi(r)}{r^2}\}$ and the overlap will have size $\leq 2\delta$ where $\delta :=\min\{ \frac{\psi(q)}{q^2},\frac{\psi(r)}{r^2}\}$.
 Now the $\frac aq - \frac br$ are in 1-to-1 correspondence with the $\frac n{qr}$ as $n$ runs through the reduced residue classes mod $qr$. Therefore, by the small sieve,
 \begin{align*}  
 \mu(E^*_q\cap E^*_r)&\leq 2\delta \#\{ n: |n|\leq 2\Delta qr \text{ and } (n,qr)=1\} \ll
 \delta\Delta qr \prod_{\substack{p|qr \\ p\leq \Delta qr}}\bigg( 1 -\frac 1p \bigg) \\
 &\leq \frac{\phi(q)\psi(q)}{q^2} \cdot   \frac{\phi(r)\psi(r)}{r^2} \cdot \exp\bigg( \sum_{\substack{p|qr \\ p> \Delta qr}} \frac 1p \bigg)  \ll   \mu(E^*_q)\mu(E^*_r) \exp\bigg( \sum_{\substack{p|qr \\ p> \Delta qr}} \frac 1p \bigg) .
 \end{align*}  
 (If $(q,r)>1$ then we need only alter this by taking $p| qr/(q,r)^2$ instead of $p|qr$ in the sum over $p$ on the far right of the previous displayed equation.)
 
 Using this one can easily  deduce the Duffin-Schaefer conjecture provided $\psi(\cdot)$ does not behave too wildly. For example Erd\H os and Vaaler \cite{Er, Vaa} proved the Duffin-Schaefer  conjecture provided the $\psi(n)$ are bounded. Key to this is to note that there are $\ll e^{-y}x$ integers $n\leq x$ for which
 \[
 \sum_{\substack{p|n \\ p> y}} \frac 1p \geq 1.
 \]
  Therefore we obtain good enough bounds on $\mu(E^*_q\cap E^*_r)$ in the previous displayed equation    
  unless $(q,r)$ is large, and unless $q$ and $r$ are each divisible by a lot of different small prime factors.  This  reduces the problem to one in the \emph{anatomy} of integers (a concept that is brought to life in the graphic novel \cite{GGL}).

\subsection{The anatomy of integers}
By partitioning $[Q,R]$ into dyadic intervals and studying the contribution of the integers in such intervals to the total we find ourselves drawn towards the following
\medskip

\noindent \textbf{Model Problem} \emph{Fix $\eta\in (0,1]$. Suppose that $S$ is a set of $\gg \eta Q/B$ integers  in $[Q,2Q]$ for which there are at least $\eta |S|^2$ pairs $q,r\in S$ such that $(q,r)\geq B$.
Must there be an integer $g\geq B$ which divides $\gg_\eta  Q/B$ elements of $S$?}
\medskip

The model problem is false but a technical variant, which takes account of the $\phi(q)/q$-weights, is true.\footnote{Let $Q=\prod_{p\leq 2y} p$  and $S:=\{ Q/p: y< p\leq 2y\}$. If $q=Q/p,r=Q/\ell\in S$ then
$(q,r)=Q/p\ell\geq B:=Q/4y^2$, but any integer $\geq B$ divides no more than two elements of $S$. (This is adapted from an idea of Sam Chow.)}
Using this one can reduce the problem to the Erd\H os-Vaaler argument, by anatomy of integers arguments, and prove the theorem.

To attack the (variant of the) Model Problem, Koukoulopoulos and Maynard view it as a question in graph theory: 
\subsection{Graph Theory}  Consider the graph $G$, with vertex set $S$ and edges between vertices representing pairs of integers with gcd$>B$.
  
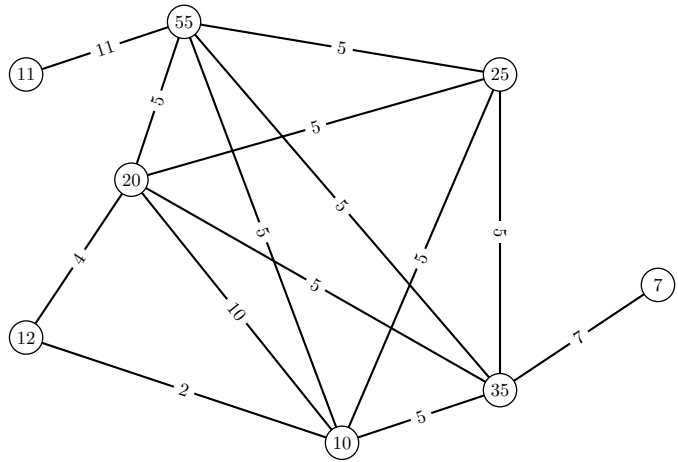
\begin{figure}[H]
\begin{center}
\begin{tikzpicture}[scale=.7,transform shape]
  \Vertex[x=0,y=7]{11}
  \Vertex[x=0,y=2]{12}
  \Vertex[x=2,y=5]{20}
  \Vertex[x=3,y=8]{55}
  \Vertex[x=6,y=0]{10}
  \Vertex[x=9,y=7]{25}
  \Vertex[x=9,y=1]{35}
  \Vertex[x=12,y=3]{7}
  \tikzstyle{LabelStyle}=[fill=white,sloped]
  \Edge[label=$11$](11)(55)
 \Edge[label=$5$](20)(55)
 \Edge[label=$5$](55)(25)
 \Edge[label=$5$](25)(35)
 \Edge[label=$7$](35)(7)
 \Edge[label=$5$](10)(35)
 \Edge[label=$2$](12)(10)
 \Edge[label=$4$](12)(20)
 \Edge[label=$10$](20)(10)
 \Edge[label=$5$](55)(10)
 \Edge[label=$5$](10)(25)
 \Edge[label=$5$](20)(35)
 \Edge[label=$5$](20)(25)
 \Edge[label=$5$](55)(35)
\end{tikzpicture}
\caption{Vertices $=$ The integers in our set. \newline
.\hskip.8in Edges $=$  Pairs of integers with a large GCD.}
\end{center}
\end{figure}

\noindent Beginning with such a graph for which the edge density is $\eta$,  we wish to prove that there is a ``dense subgraph'' $H$ whose  vertices are each divisible by a fixed integer $\geq B$. To locate this structured subgraph $H$,  Koukoulopoulos and Maynard use an 
iterative ``compression'' argument,  inspired by the papers of Erd\"os-Ko-Rado \cite{EKR} and Dyson \cite{Dys}: with each iteration, they pass to a smaller graph but with more information about which primes divide the vertices. This is all complicated by the weights $\phi(q)/q$. The details are complicated (see a vague sketch in the next subsection); and the reader is referred to \cite{Kou}, where the original proof of \cite{KM} is better understood from more recent explorations of Green and Walker \cite{GW}, who gave an elegant proof of the following important variant:

\emph{If $R\subset [X,2X]$ and $S\subset [Y,2Y]$ are sets of integers for which $(r,s)\geq B$ for at least 
$\delta |R||S|$ pairs $(r,s)\in R\times S$ then $|R||S|\ll_\epsilon \delta^{-2-\epsilon} XY/B^2$.}

Although this has a slightly different focus from the model problem, it focuses on the key question of how large such sets can get and takes account of the example of footnote 8 (unlike the model problem).

\subsection{Iteration and graph weights} The key to such an iteration argument is to develop a measure which exhibits  how close one is getting to the final goal, which can  require substantial ingenuity. In their paper
Koukoulopoulos and Maynard \cite{KM} begin with two copies of $S$ and construct a bipartite graph $V_0\times W_0$ with edges in-between $q\in V_0=S$ and $r\in W_0=S$ if $(q,r)\geq B$. The idea is to select distinct primes $p_1, p_2,\ldots$ and then $V_j=\{ v\in V_{j-1}: p_j \text{ divides } v\}$ or 
$V_j= \{ v\in V_{j-1}: p_j \text{ does not divide } v\}$, and similarly $W_j$,  so that $p_j$ divides all $(v_j,w_j), v_j\in V_j, w_j\in W_j$ or none. If we terminate at step $J$ then there are integers $a_J,b_J$, constructed out of the $p_j$, such that $a_J$ divides every element of $V_J$ and  $b_J$ divides every element of $W_J$. The goal is to proceed so that $(v_J,w_J)\geq B$ for some $J$, for all
$v_J\in V_J, w_J\in W_J$ such that  all of the prime divisors of any $(v_J,w_J)$ appears amongst the $p_j$. Hence, if say all the integers in $S$ are squarefree, then $(a_J,b_J)=(v_J,w_J)\geq B$. So how do we measure progress in this algorithm?

One key measure is $\delta_j$, the proportion of pairs $v_j\in V_j, w_j\in W_j$ with $(v_j,w_j)\geq B$, another the size of the sets $V_j$ and $W_j$. Finally we want to measure how much of the $a_jb_j$ are given by prime divisors not dividing $(a_j,b_j)$, which we can measure using $\frac{a_jb_j}{(a_j,b_j)^2}$. Koukoulopoulos and Maynard \cite{KM} found, after some trial and error, that the measure
\[
\delta_j^{10} \cdot |V_j|  \cdot |W_j|  \cdot \frac{a_jb_j}{(a_j,b_j)^2}
\]
fits their needs, allowing them eventually to restrict their attention to $v,w\in S$ for which $a_j$ divides $v$, $b_J$ divides $w$ and 
\[
 \sum_{\substack{p| vw/(v,w)^2 \\ p> y}} \frac 1p \approx 1.
 \]
Koukoulopoulos and Maynard then finish the proof by applying a relative version of the Erd\" os-Vaaler argument to the pairs $(v/a_J, w/b_J)$.

\subsection{Hausdorff dimension} If $\sum_{n\geq 1}    \phi(n)  \cdot  (\psi(n) /n^2)$ is convergent then $\mu(\mathcal L^*(\psi))=0$ so we would like to get some idea of the true size of  $\mathcal L^*(\psi)$.
Using a result of Beresnevich and Velani \cite{BV}, one can deduce that   the Hausdorff dimension of $\mathcal L^*(\psi)$ is given by the infimum of the real $\beta> 0$ for which 
\[
\sum_{n\geq 1}    \phi(n)  \cdot  \bigg( \frac{ \psi(n) }{n^2}   \bigg)^\beta
\text{  is convergent.}
\]


\bibliographystyle{amsplain}

\end{document}